\documentclass[11pt]{article}
\usepackage{geometry}                
\usepackage{graphicx}
\usepackage{amsmath}
\usepackage{amsthm}
\usepackage{amssymb}
\usepackage{leftidx}
\usepackage{epstopdf}
\usepackage{comment}

\newtheorem{theorem}{Theorem}

\begin{document}

\newcommand {\dd}{\textup{d}}
\newcommand{\RR}{\mathbb{R}}
\newcommand{\slipparameter}{\beta}
\newcommand{\eigparam}{\lambda}
\newcommand{\remaind}{{\rm rem}}
\newcommand{\steady}{{\rm steady}}

\newcommand{\CM}{{{\mathbf C}\mathbf M}}
\newcommand{\AM}{{{\mathbf A}\mathbf M}}

\numberwithin{equation}{section} 

\title{Some isoperimetric results concerning undirectional flows in microchannels}
\author{Grant Keady and Benchawan Wiwatanapataphee\\
Department of Mathematics and Statistics\\
Curtin University, Bentley, 6102\\
Western Australia}
\date{\today}                                           

\maketitle
\vspace{0.3cm}
\section*{Abstract}
Three isoperimetric results are treated. 
(i) At a given pressure gradient, for all channels with given (cross-sectional) area that which
maximises the steady flow $Q_\steady$ has a circular cross-section.
(ii) Consider flows starting from prescribed initial conditions developing from a prescribed imposed pressure gradient,
either periodic or steady.
For such flows, amongst all channels with given  area, that which generically
has the slowest approach to the long-term, periodic or steady, flow is the circular disk cross-section.
(iii) Similar results for polygonal, $n$-gon, channels, with the optimising shape being the regular $n$-gon are discussed.\\
This arXiv preprint will supplement the journal paper (submitted just before this supplement):
the journal paper reports, concisely, the isoperimetric results of Theorems 1, 2 and 3.
Items here additional to the journal paper include further isoperimetric results,
estimates involving geometric functionals besides area such as perimeter and moment of inertia,
perturbation analysis of nearly circular domains, and reporting on some previously published conjectures.


\section{Introduction}
\label{sec:Intro}

\subsection{Overview}
\label{subsec:Overview}

General remarks on microchannels and applications are given in many papers,
in particular 
\cite{KWW13,SWLW13,SWLW14,WWS14}.
The papers just cited develop explicit exact solutions of pde problems,
similar to those of \S\ref{subsec:pde}, for various simple geometries for the cross-sectional shape
$\Omega$ of the microchannel;
rectangles in~\cite{WWS14},
circular cross-section in~\cite{SWLW13}.
This paper has a different focus, namely that of describing qualitative properties and bounds
on solutions to functionals of solutions of the pde problem.
The bounds in the main theorems of this paper are universally known as `isoperimetric' inequalities, though we remark that
here they are actually `iso-area'.
Both explicit solutions and bounds are useful for checking the results from numeric codes
used by engineers.

The new announcements in this paper concern slip boundary conditions,
equation~(\ref{eq:slipbc}) with $\beta>0$.
The no-slip case has $\beta=0$.
The results stated in the abstract are well-known, from around 1950,
for no-slip boundary conditions.
See ~\cite{PoS51}.
Although some engineers have {\it assumed} --
correctly as it turned out --
that they are also true with slip boundary conditions,
the proof that this is the case is relatively recent.
The papers containing the proofs of Theorems~\ref{thm:thm1} and~\ref{thm:thm2} have other applications in mind --
heat flow with Newton's law of cooling, for example --
and the main new element in this paper is to call attention to the results in connection with microchannels.
Theorem~\ref{thm:thm3} -- a further observation concerning rectangular microchannels -- is new, 
but its proof is merely calculation, building from calculations similar to those in~\cite{WWS14}.
The result  of Theorem~\ref{thm:thm3} is
important mostly as evidence that some open problems already stated in the literature
may be answered affirmatively.

The proofs for the no-slip case use symmetrisation: see~\cite{PoS51}.
However symmetrisation is not appropriate for the slip-flow case.
Partly to reassure engineer readers who may find the proofs of the theorems difficult,
but who know the explicit exact solutions, we use these to illustrate the main results.
Those readers who just want the main results  need only read
\S\ref{sec:Intro}, \S\ref{subsec:isoperimetricQSteady} and \S\ref{subsec:isoperimetricLambda1}
and can omit all the explicit solutions.
Some of the explicit solution work is perhaps helpful in connection with some of the details
in the discussion and open questions of   \S\ref{sec:Discussion} and \S\ref{subsubsec:PolygonalbGt0}.

\subsection{The pde problem}
\label{subsec:pde}

Batchelor's classic book on fluid mechanics~\cite{Ba67} \S4.3 
gives a number of solutions for unsteady unidirectional flows of a viscous fluid
with no-slip boundary conditions.
See also~\cite{Wa89}.
(Steady solutions with no slip are treated in~\cite{Ba67} \S4.2 
and many papers, e.g.~\cite{TB11}, 
and in the context of microchannels with slip
in~\cite{BYC,BTT09,DM07,LBS05,TH11,Wa12,Wa12b}.)
We use coordinates in which $z$ is in the direction of the flow, and $x$ and $y$ transverse to it,
$(x,y)$ varying over a cross-section $\Omega$.
The pressure gradient $p_z(t)$ is prescribed.
Then equations (4.2.2) and (4.2.3) from~\cite{Ba67} become
\begin{equation}
\rho \frac{\partial u}{\partial t}
= {p_z}(t) + \mu \left( \frac{\partial^2 u}{\partial x^2}+ \frac{\partial^2 u}{\partial y^2}\right) 
\quad {\mbox{\rm for }}\ \ (x,y)\in\Omega\ \  {\mbox{\rm and }}\ \  t>0 .
\label{eq:pde}
\end{equation}
The notation is as in~\cite{Ba67}:
$u$ is the fluid velocity, $\rho$ density and $\mu$ viscosity.
The boundary condition we treat is, with $\beta\ge{0}$,
\begin{equation}
u +\slipparameter \frac{\partial u}{\partial n}
= 0 \qquad {\mbox{\rm on }}\ \ \partial\Omega  ,
\label{eq:slipbc}
\end{equation}
where $n$ denotes the outward normal to the boundary.
In the context of flows 
this boundary condition is called a  {\it slip boundary condition} or
{\it Navier's boundary condition}; 
in the wider mathematical literature it is called a {\it Robin boundary condition}.
In~\cite{Ba67} it is just the $\slipparameter=0$ case that is studied.

It is appropriate here to note that how well the Navier (Robin) boundary condition
fits real microchannel flows is problematic,
see~\cite{NW11}.
It is an open question as to how one might generalize to obtain
results similar to those of Theorems 1 and 2
when boundary conditions more accurately modelling real flows are used.

Absorb the $\rho$ and $\mu$ into  space and time variables so that we take $\rho=1$ and $\mu=1$.

There are several functionals of interest, for example the volume flux
\begin{equation}
 Q(\slipparameter,t)= \int_\Omega u .
\label{eq:Qdef}
\end{equation}

Of course the same pde problem arises in other contexts.
One of these is heat flow where the Robin boundary condition is known as Newton's law of cooling.

The special cases of ${p_z}(t)$ we study here are\\
$\bullet$ steady flows, where $p_z$  is positive and constant in time, and $u$ is also constant in time;\\
$\bullet$ periodic flows (which can be found by superposition of solutions) when $p_z=\exp(i\omega t)$  and $u(x,y,t)=u_p(x,y)\exp(i\omega t)$ is also periodic with the same period (and the steady flow case is that of $\omega=0$);\\
$\bullet$ transient flows where $p_z$ is a Heaviside step function in time multiplying the functions above and initial data $u(x,y,0)$ is also prescribed.
When $p_z$ is a Heaviside step function,
 $p_z(t)={\Delta{p}}\ {\rm Heaviside}(t)$ with ${\Delta{p}}>0$, and  $u(x,y,0)=0$  we refer to this as a `starting flow'.
For a starting flow, and the details of the initial data do not matter, $Q_\steady(\slipparameter)=Q(\slipparameter,\infty)$.

What is new in the main part of this paper is not the Theorems~\ref{thm:thm1} and~\ref{thm:thm2} stated here --
these were proved by others --
but noting their application to microchannel flows with slip.
The same pde problems, but in polygonal domains, are treated in the appendices:
Theorem~\ref{thm:thm3} concerns transient flows in rectangular microchannels.

\subsection{Eigenfunction expansions}
\label{subsec:eigfn}

The linear constant coefficient  pde problems here 
are standard undergraduate applied mathematical exercises, the methods being treated in
many text books, e.g.
Courant and Hilbert~\cite{CHI57}, 
Lebedev at al.~\cite{LSU79}, etc.

Sturm-Liouville style eigenvalue problems arise on seeking solutions of the heat 
equation 
 in the form
$\exp(-{\lambda_n}t)\phi_n(x,y)$.
This leads to the Helmholtz eigenvalue problem to find the
$(\lambda_n,\phi_n)$, with $\phi_n$ not identically zero, satisfying
\begin{equation}
\left( \frac{\partial^2 \phi_n}{\partial x^2}+ \frac{\partial^2 \phi_n}{\partial y^2}\right) 
+\lambda_n \phi_n = 0.
\label{eq:Helmholtz}
\end{equation}
The $\lambda_n$ are ordered so that $\lambda_{n+1}\ge{\lambda_n}$.
Completeness of the eigenfunctions and expansions in terms of them is
treated in classical mathematical methods books such as~\cite{CHI57}.

\subsubsection{Periodic and steady flows}

Suppose ${p_z}(t)=\exp(i \omega t)$.
Then seek solutions in the form
$$ u(x,y;t)
=  \sum_{j=1}^\infty C_j \exp(i \omega t) \phi_j(x,y) , $$
where here $\phi_j$ is normalized in the usual way.
On substituting this into equation~(\ref{eq:pde}) one gets
$$
i \omega C_j = \langle 1,\phi_j\rangle - \lambda_j C_j $$
so that
$$ C_j = \frac{\int_\Omega \phi_j}{\lambda_j + i\omega} .$$
Hence
\begin{equation}
Q(\slipparameter,t)
=\left(\sum_{j=1}^\infty
\frac{1}{\lambda_j + i\omega}
 \frac{(\int_\Omega \phi_j)^2}{\int_\Omega \phi_j^2}.\right)\exp(i\omega t) .
\label{eq:QexpPeriodic}
\end{equation}

On letting $\omega$ tend to zero, we get the steady-flow case
\begin{equation}
\frac{u_\steady(x,y) }{{\Delta{p}}}
=  \sum_{j=1}^\infty \frac{\int_\Omega \phi_j}{\lambda_j\int_\Omega \phi_j^2} \phi_j(x,y) \ .
\label{eq:uSteadySeries}
\end{equation}
Integrating equation~(\ref{eq:uSteadySeries}) or letting  $\omega$ tend to zero
in equation~(\ref{eq:QexpPeriodic}) gives, with
\begin{equation}
Q_j
=\frac{( \int_\Omega \phi_j)^2 }{\lambda_j \int_\Omega \phi_j^2} \  ,
\label{eq:Qjdef}
\end{equation}
that
\begin{equation}
\frac{Q_\steady(\slipparameter) }{{\Delta{p}}}
=  \sum_{j=1}^\infty Q_j \ .
\label{eq:QQjSeries}
\end{equation}

There are methods other than this eigenfunction method to solve for periodic and for steady flows.
These lead to rather different formulae, as exemplified in the examples of $\Omega$
circular or rectangular or equilateral triangular given later in this paper.
It then becomes an exercise to demonstrate that these formulae agree with the
eigenfunction expansions.

\subsubsection{Transients}

We now treat initial boundary-value problems. 
After subtracting off the flow that would be achieved at large time, the time-dependent remainder $u_\remaind(x,y,t)$
satisfies the heat equation
\begin{equation}
\frac{\partial u_\remaind}{\partial t}
=  \left( \frac{\partial^2 u_\remaind}{\partial x^2}+ \frac{\partial^2 u_\remaind}{\partial y^2}\right) 
\quad {\mbox{\rm for }}\ \ (x,y)\in\Omega\ \  {\mbox{\rm and }}\ \  t>0 ,
\label{eq:heateq}
\end{equation}
together with the slip boundary condition~(\ref{eq:slipbc}) and given initial values.

This leads to developing the solution
 $u_\remaind$
as an eigenfunction expansion
\begin{equation}
{u_\remaind}(x,y,t)
=\sum_{j=1}^\infty \frac{\int_\Omega {u_\remaind}(x,y,0)\phi_j}{\lambda_j\int_\Omega \phi_j^2} \phi_j(x,y) \exp(-\lambda_j t))\ .
\label{eq:urGeneral}
\end{equation}

One particular example is that of `starting flows' where the step-function pressure gradient is turned on at $t=0$
acting on an initially still liquid:
\begin{equation}
u(x,y,t)
= u_\steady(x,y) - {{\Delta{p}}} \sum_{j=1}^\infty 
\frac{\int_\Omega  \phi_j}{\lambda_j\int_\Omega \phi_j^2} 
 \exp(-{\lambda_j}t)\phi_j(x,y) \ .
\label{eq:uStartingSeries}
\end{equation}
From equation~(\ref{eq:uSteadySeries}) one finds
\begin{equation}
Q(\slipparameter,t) 
= Q(\slipparameter,\infty) -{\Delta{p}}  \sum_{j=1}^\infty \exp(-{\lambda_j}t) 
Q_j .
\label{eq:QphiSeries}
\end{equation}
(We remark that the normalisation of $\phi$ is not needed in equations~(\ref{eq:uStartingSeries}-\ref{eq:QphiSeries}).)
The function $Q(\slipparameter,\cdot)$ is positive, increasing, concave for $t\in[0,\infty)$.\\
For large time we have
\begin{equation}
Q(\slipparameter,t) 
\sim Q(\slipparameter,\infty) -  {\Delta{p}} \exp(-{\lambda_1}t) Q_1\ 
{\rm as\ } t\rightarrow\infty\ ,
\label{eq:Qlarget1}
\end{equation}
with $Q_1$ being as in equation~(\ref{eq:Qjdef}) with $j=1$.
The series~(\ref{eq:QQjSeries})  for $Q_\steady(\beta)$ is evident from
equation~(\ref{eq:QphiSeries}) on setting $t=0$ and
noting $Q(0,\beta)=0$:
See also~\cite{PoS51}p106 equation (12) (and for $\Omega$ a circle, p109, equation (21)).
Some analytic solutions for starting flows in simple geometric shapes are given
in~\cite{Ng16}.

More generally, the large-time asymptotics is dominated by the periodic or steady flow at large time and, generically, the term in
$\exp(-\lambda_1 t)$,
 $\lambda_1$ determines the rate at which the large time solution is approached.
 The generic qualification is because in very exceptional circumstances the coefficient of
 $\exp(-\lambda_1 t)$ in equation~(\ref{eq:urGeneral}) might be zero.

Some analytic solutions for starting flows in simple geometric shapes are given
in~\cite{Ng16}.

\subsection{Overview, ctd.}\label{subsec:OverviewCtd}

Our main results are Theorems 1 to 3.
\begin{itemize}
\item Theorem~\ref{thm:thm1} is in \S\ref{sec:Steady},
Theorem~\ref{thm:thm2} in \S\ref{sec:Transient}.\\
We remark that all the sections of the main part from \S\ref{sec:Symmetrisation}
to \S\ref{sec:nearCirc} record matters related to attempts to extend the results in various directions.
\item
Also of note, though concerning just $\lambda_1$ for rectangles, is
Theorem~\ref{thm:thm3} in Appendix~\ref{subsec:Proofexplicit}.\\
As this is a new result, we have attempted different proofs of it, and also allowed ourselves to
establish other inequalities on $\lambda_1$.
\end{itemize}

The table below provides a short and partial `index'  of these notes, intended for reference after some study of the notes.
\medskip

\begin{tabular}{|| c | c | c | c | c | c ||}
\hline
\S& geometry & Comments& $Q(t)$ & $Q_\steady$& $\lambda_1$ \\
\hline
& & & & & \\
\S\ref{sec:circular}& Circle& & &\S\ref{subsec:circUltimate} & \S\ref{subsec:circStarting}\\ 
& & & & & \\
\S\ref{sec:polygonal}& EquilateralTriangle& & &\S\ref{subsec:equilatQs}&\S\ref{subsec:equilatLambda1}\\
\S\ref{sec:polygonal}& Rectangle& & &\S\ref{subsec:rectQs}&\S\ref{subsec:Rectangular}\\
& & & & & \\
\S\ref{sec:Steady}& &Theorem 1 & & * & \\
\S\ref{sec:Transient}& &Theorem 2 & & & * \\
& & & & & \\
\S\ref{sec:IllustrativeTimeDep}& Circle& & * & & \\
\hline
\end{tabular}
\medskip

As noted above,  the sections of the main part from \S\ref{sec:Symmetrisation}
to \S\ref{sec:nearCirc} record matters related to attempts to extend the results in various directions.
Specifically \S\ref{sec:Symmetrisation} notes that symmetrisation is appropriate when $\beta=0$, but
inappropriate when $\beta>0$.
After this \S\ref{sec:Perturbations} introduces the following two sections and notes that
improvements to some isoperimetric inequalities can be effected using moments of inertia about the centroid.
The next section~\S\ref{sec:Measuring} addresses the question of how the departure of the domain functional,
especially $Q_\steady(\Omega)$ and $\lambda_1(\Omega)$ (but also simpler functionals such as
polar moments of inertia about the centroid) might be estimated in terms of simple geometric quantities --
modulus of asymmetry and isoperimetric deficit being examples.
We use the ellipse as an example.
In~\S\ref{sec:nearCirc} we treat domains which are nearly circular.
Again, the main example we treat is the ellipse.
There is a discussion and mention of some open questions in~\S\ref{sec:Discussion}.

The appendices treat polygonal domains.
Open questions concerning general polygons, $n$-gons, some due to Polya and Szego~\cite{PoS51} for $\beta=0$,
are described in Appendix~\ref{sec:PolygonalIso}.
When $\beta=0$ \cite{PoS51} have results for $n=3$, triangles, and $n=4$ quadrilaterals.
A small result towards $n=4$ and $\beta>0$, for rectangles is given in our Theorem~\ref{thm:thm3} proved in
Appendix~\ref{sec:rectlambda1}.
The final appendix~\ref{sec:compmon} treats some technicalities arising from side-issues to Theorem~\ref{thm:thm3}

\section{Circular channels}\label{sec:circular}

We don't actually need the explicit formulae of this section (or the next)
for the theorems below,  in~\S\ref{subsec:isoperimetricQSteady} and in~\S\ref{subsec:isoperimetricLambda1}.
We offer the following partial justification.
(i) Our theorems compare flows in different cross-sections to those in circular channels.
(ii) The formulae~(\ref{eq:circuStarting},\ref{eq:circQStarting}) below do provide an example of the eigenfunction expansions
used in this paper.
(iii) Ratios of domain functionals ($Q_\steady$ and $\lambda_1$) comparing other shapes 
(e.g. the rectangular shapes of ~\cite{WWS14})
with those of the circle with the same area are given in
\S\ref{sec:polygonal}
and provide some indication that the theorems with $\beta>0$ are more delicate than the
classical, 1950s, results at $\beta=0$.

The isoperimetric inequalities of this paper are in connection with
some functional being extreme when the problem is solved in a circular domain.
A classic reference on isoperimetric inequalities is~\cite{PoS51}
and the results in that book establish the theorems mentioned in this paper
when $\beta=0$.

\subsection{Circular cross-section: ultimate behaviour}\label{subsec:circUltimate}

When $\Omega$ is a circular disk of radius $a$, the solutions are readily found using polar coordinates $r=\sqrt{x^2+y^2}$.
The periodic flow is
$u(r,t)=u(r)\exp(i\omega t)$ where $u(r)$ satisfies the o.d.e.
$$ i\omega u
= 1 +\frac{1}{r}\frac{\partial}{\partial r} r\frac{\partial u}{\partial r} .
$$
Define 
$$ \sigma=\frac{\left( 1-\,i \right)}{\sqrt{2}} \sqrt {\omega} \qquad
{\mbox{\rm so\ that\ }}
\sigma^2=-i\,\omega, \qquad \omega=i\sigma^2 .
$$
Using  the boundary conditions that $u(0)$ is bounded, this and the o.d.e. are satisfied by
\begin{equation}
u(r)
=\frac{ {C J_0\left(\sigma r\right)}
-{J_0\left(\sigma a\right)} }
{  {J_0\left(\sigma a\right)} \sigma^2}
= \frac{ {C J_0\left(\sigma r\right)} }
{  {J_0\left(\sigma a\right)} \sigma^2} -\frac{1}{\sigma^2} .
\label{eq:oscCircGen}
\end{equation}

There is less clutter in the no-slip case, i.e. when $\beta=0$, is
The solution of the preceding ode using the boundary conditions $u(0)$ bounded and $u(a)=0$ is
obtained by setting $C=C_0=1$ in equation~(\ref{eq:oscCircGen}).
The oscillatory volume flux is $Q(t)=q\exp(i\omega t)$ where
\begin{eqnarray*} q
&=& 2\pi\int_0^a u(r) r \, dr \\
&=& -{\frac {\pi \,a\, \left( -2\,{J_1\left(\sigma\,a\right)}+a\,\sigma\,{J_0\left(\sigma\,a\right)} \right) }
{{\sigma}^{3}{J_0\left(\sigma\,a\right)}}} 
\quad\mbox{\rm when\ } \beta=0\ .
\label{eq:qOscbEq0}
\end{eqnarray*}

When $\beta>0$ the solution is given by equation~(\ref{eq:oscCircGen})
where now 
$$ C = C_\beta
=  \frac{J_0(\sigma{a})}{J_0(\sigma{a}) -\beta\sigma J_1(\sigma{a})}.
$$
Now the volume flux parameter $q$ is
\begin{equation}
q
={\frac {\pi \,a\, \left( -2\,{J_1\left(\sigma\,a\right)}+a\,\sigma\,{J_0\left(\sigma\,a\right)} -
\sigma^2 a\beta J_1\left(\sigma a\right)\right) }
{{\sigma}^{3}
({J_0\left(\sigma\,a\right)}-\beta\sigma J_1\left(\sigma\,a\right))}}.
\label{eq:qOscbNe0}
\end{equation}

The steady flow is most simply derived directly, but it can also be obtained as the limit of the preceding
periodic solutions when $\omega$ tends to zero. 
The steady flow is
\begin{equation}
 \frac{u_\steady}{{\Delta{p}}} 
= \frac{1}{4}\left(a^2-r^2\right) +\frac{\slipparameter a}{2} .
\label{eq:circuSteady}
\end{equation}
When $\slipparameter=0$ this is Poiseuille flow.

For general $\beta\ge{0}$ we have
\begin{equation}
Q_\steady(\slipparameter)=Q(\slipparameter,\infty)
=  \frac{\pi a^3}{8} (a+4\slipparameter) {\Delta{p}} . 
\label{eq:circleQinf}
\end{equation}
See also~\cite{La32}\S{331},~p586.
Equation~(\ref{eq:circleQinf}) is very well-known and is used
in the experimental determination of $\beta$.

\subsection{Circular Cross-Section: Transients}\label{subsec:circStarting}

The eigenvalues corresponding to radial eigenfunctions are found by separation of variables.
These eigenvalues are $\lambda_j=(\gamma_j/a)^2$  where the $\gamma_j$ are consecutive positive roots of
\begin{equation}
J_0(\gamma)=\frac{\slipparameter}{a} \gamma J_1(\gamma) . 
\label{eq:JrootsBetaNe0}
\end{equation}
It is the least positive solution $\gamma_1$ , which dominates at large time:
with 
\begin{equation}
\lambda_1 = \left(\frac{\gamma_1}{a}\right)^2 .
\label{eq:lambda1gamma}
\end{equation}

\begin{figure}[hb]
\centerline{\includegraphics[height=7cm,width=13cm]{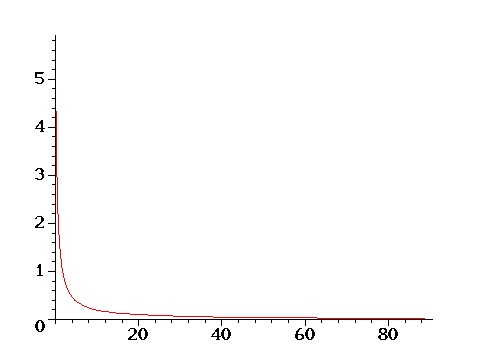}}
\caption{Plot of $a^2\lambda_1=\gamma_1^2$ against  $\slipparameter/a$.}
\label{fig:lambdabetaplot}
\end{figure}

A plot showing the behaviour of $\lambda_1$ as $\beta$ varies is given in
Figure~\ref{fig:lambdabetaplot}.
When $\beta=0$, $\lambda_1=(j_0/a)^2$ and as $\beta$ tends to infinity, $\lambda_1$ tends to zero.
We have
$$\left(\frac{\gamma_1}{a}\right)^2 =\lambda_1
\sim \frac{|\partial\Omega|}{\beta|\Omega|} = \frac{2}{\beta a}\qquad
{\rm as\ }\ \beta\rightarrow\infty .
$$

\clearpage

As a simple example of flows involving transients, consider starting flows, with the fluid initially at rest.
The series solution, when $\slipparameter=0$ is given in Batchelor~\cite{Ba67} \S4.3 equation~(4.3.19).
The series solution when $\slipparameter\ge{0}$ is given (in the context of heat flow) in Problem 247 of~\cite{LSU79}, 
and results there include the following:
\begin{equation}
\frac{u(r,t)}{{2a^2 \Delta{p}}}
= \frac{1}{8}(1-\frac{r^2}{a^2} ) +
\frac{\slipparameter}{4 a} -
\sum_{j=1}^\infty \frac{J_0(\gamma_j \frac{r}{a})}
{\gamma_j^2(1+(\frac{\slipparameter}{a}\gamma_j)^2 )J_1(\gamma_j)}
\exp(-\gamma_j^2 \frac{t}{a^2})
 .
\label{eq:circuStarting}
\end{equation}

Integrating~(\ref{eq:circuStarting}) gives

\begin{equation}
\int_0^a r J_0(\gamma r)\, dr
=\frac{a^2\, J_1(\gamma)}{\gamma} ,
\label{eq:besselInt}
\end{equation}
and hence
\begin{equation}
\frac{Q}{{a^2 \Delta{p}}}
= \frac{\pi a}{8}(a+4\slipparameter)-
4\pi a^2 \sum_{j=1}^\infty \frac{\exp(- \gamma_j^2 t/a^2)}
{\gamma_j^3 (1 + (\frac{\slipparameter}{a}\gamma_j)^2) }
 .
\label{eq:circQStarting}
\end{equation}



\section{Triangular and rectangular cross-sections}\label{sec:polygonal}

The formulae in this section are used to compare their domain functionals with those of a circle with the same area.
(The results  are briefly reviewed in \S\ref{subsec:QscircleTriSq}
and in
\S\ref{subsec:lambda1circleTriSq}
in connection with illustrating theorems.
They are also referenced in \S\ref{subsec:QtcircleTSq}.)
They are also used in
Appendix~\ref{sec:PolygonalIso}
where we consider triangles with a given area and
rectangles with a given area.


\subsection{$Q_\steady$ for equilateral triangles}\label{subsec:equilatQs}

There is a simple polynomial solution for the steady flow when $\Omega$ is an equilateral triangle.
We suspect that it may date from the nineteenth century, but, 
unlike the situation for $\lambda_1$ we do not know papers treating it.
 We treat the triangle with vertices $(-a,0), (a,0), (0,a\sqrt{3})$.
The elastic torsion $\slipparameter=0$ solution has been known for centuries:
$$\frac{u_{\Delta,0 }} {\Delta{p}}
= \frac{y \left( y -\sqrt{3}( a +x) \right) \left( y -\sqrt{3}( a -x )\right)}{4 \sqrt{3} a}\qquad
{\mbox{\rm for }}\ \slipparameter=0 .
$$
See equation~(3.12) of~\cite{MK94}. Define
$$\frac{u_{\Delta,\infty}} {\Delta{p}}
= \left(\frac{a^2}{6}-\frac{1}{4}
\left( x^2 + (y-\frac{a}{\sqrt{3}})^2  \right)\right) .$$
Then, 
as in equation~(3.16) of~\cite{MK94},
$$\frac{u_{\rm steady}}{\Delta{p}}
=\frac{u_{\Delta,\slipparameter }}{\Delta{p}}
= \frac{ \frac{a}{\sqrt{3}} u_{\Delta,0} + \slipparameter u_{\Delta,\infty}+\frac{a}{2\sqrt{3}}  \slipparameter (\frac{a}{\sqrt{3}} + \slipparameter)}
{\frac{a}{\sqrt{3}} + \slipparameter} .
$$
The solution has been found, independently, many times and reported frequently: e.g.~\cite{Wa03,Wa12}.

\medskip
Formulae for  $Q_\steady$ have been found several times.
When $\beta=0$ it is a classical formula for the torsional rigidity of an equilateral triangle, due to St Venant,
and we denote it here, as in~\cite{MK94}  by $S_0$, and
$$ S_0 = \frac{\sqrt{3} a^4}{20} = \frac{|\Omega|^2}{20\sqrt{3}} .$$
The formula for $Q_\steady$ is conveniently written in terms of $S_0$ and
$S_\infty$ where $S_\infty$ is the integral of $u_\infty$ over $\Omega$:
$$ S_\infty = \frac{a^4}{4\sqrt{3}} =  \frac{|\Omega|^2}{12\sqrt{3}} .$$
Finally, with the area $|\Omega|= a^2\sqrt{3}$ and 
perimeter $|\partial\Omega|= 6a$,
\begin{equation}
 \frac{Q_\steady}{\Delta{p}} 
 = \frac{\beta|\Omega|^2}{|\partial\Omega|} + S_0 +
\frac{\beta}{c_0+\beta}\left(S_\infty -S_0\right) ,\qquad
\mbox{\rm where}\ c_0=\frac{a}{\sqrt{3}} = 2\frac{|\Omega|}{|\partial\Omega|}  .
\label{eq:QsEquilatTri}
\end{equation}
See equations (3.19) and (3.15) of ~\cite{MK94}.

Denote by $Q_{\steady\odot}$ the steady volume flow for the circle as given in 
equation~(\ref{eq:circleQinf}).
Denote by $Q_{\steady\Delta}$ the steady volume flow for the equilateral triangle  as given in 
equation~(\ref{eq:QsEquilatTri}).
Now set ${\Delta{p}}$ to be the same for both and als
the areas of both to be the same, and we usually take the circle to have unit radius.
One can readily plot, as a function of $\beta$,  (the rational function of $\beta$ giving) the ratio 
 $Q_{\steady\Delta}/Q_{\steady\odot}$ when both shapes have area $\pi$.
 It is of course a consequence of the St Venant inequality that at $\beta=0$ the ratio is less than one.
 It happens that when $\beta>0$ the ratio is always greater than that which it is at $\beta=0$.
 This can be established by routine algebra.
 Another simple piece of algebra involves noting that the rational function of $\beta$
 $(1-Q_{\steady\Delta}/Q_{\steady\odot})$ has two linear factors in its numerator, both of which are
 positive for $\beta\ge{0}$ and a numerator which is quadratic in $\beta$ and positive for all real $\beta$.
 This accords, of course, with the recent result given in Theorem 1 below, that the ratio 
 $Q_{\steady\Delta}/Q_{\steady\odot}$ remains  less than one for all $\beta$. 

\begin{figure}[hb]
\centerline{\includegraphics[height=5cm,width=6.5cm]{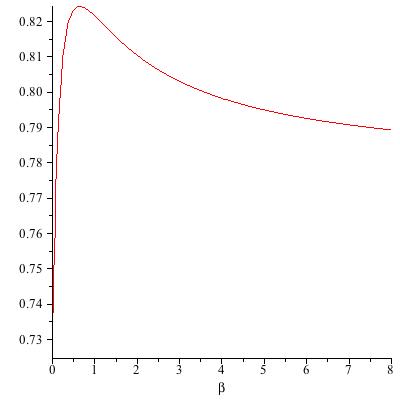}}
\caption{For domains with area $\pi$, the ratio $Q_{\steady\Delta}/Q_{\steady\odot}$
plotted against $\beta$.
For domains both with the same area, but different from $\pi$, the curves are different
but remain less than 1}
\label{fig:equilatTriQratio}
\end{figure}

\clearpage

\subsection{$\lambda_1$ for equilateral triangles}\label{subsec:equilatLambda1}

When $\beta=0$, for an equilateral triangle of side $2a$, 
from~\cite{PoS51}~p256, $\lambda_{1\Delta}=4\pi^2/(3 a^2)$.
Thus
$$\lambda_{1\Delta}=\frac{4\pi^2}{\sqrt{3}|\Omega|}
\approx\frac{22.8}{|\Omega|}
\qquad\mbox{\rm while\ }
 \lambda_{1\odot}=\frac{\pi j_0^2}{|\Omega|}
 \approx\frac{18.2}{|\Omega|}
$$
We remark that  Lam{\'e} in 1833
discovered the Dirichlet eigenvalues and eigenfunctions of the equilateral triangle,
not just the fundamental which is the main concern in this paper.

We now consider $\beta>0$.
This problem was also solved by Lam{\'e} in 1833, though the numerical study waited until better computing became available.
The history is mentioned in~\cite{McC04}.
We begin by reporting the solution in~\cite{MK94} where the equilateral triangle treated
has vertices $(-1,0), (1,0), (0,\sqrt{3})$.
The function
$$v=\sin(a+b\sqrt{3}x + b y)+\sin(a -b\sqrt{3}x + b y) +\sin(a+b\sqrt{3}-2 b y) $$
satisfies the Helmholtz equation with $\lambda=4 b^2$.
It also satisfies the Robin boundary condition on the 3 sides when the transcendental equation
\begin{equation}
 \tan\left(\frac{t\sqrt{3}}{\beta}\right)
= \frac{3 t}{2 t^2-1}\ \  {\rm and\ }
a=\arctan(t),\quad b=\frac{t}{\beta}.
\label{eq:transct}
\end{equation}
See (3.23) and (3.24) of~\cite{MK94}.
We comment that the function satisfies the pde over the whole plane, and
the Robin boundary conditions on the 3 (infinite) lines through the vertices.
The critical points on the triangle medians and line extensions of them can be calculated.
Thus
$$ v_y(0,y) = 4 b \sin\left(\frac{3 b}{2}(\frac{1}{\sqrt{3}}-y)\right) \sin\left(a +\frac{b}{2}(\sqrt{3}-y)\right) .
$$
The centroid, in-centre, of the triangle, $y=1/\sqrt{3}$, is a critical point.
More generally, $y=1/\sqrt{3} + 2 n\pi/(3 b)$ for integer $n$ are critical points.
There are other critical points at $y=(2 a + b\sqrt{3}-2 n\pi)/b$ for integer $n$.

To find $\lambda_1 (=\lambda_{1\Delta})$ we need the smallest root $t_1$ of the transcendental equation~(\ref{eq:transct})
and then 
\begin{equation}
\lambda_1=(\frac{2 t_1}{\beta})^2 .
\label{eq:lambda1TriMK}
\end{equation}
Both the left and right-hand sides of the transcendental equation are odd functions of $t$ so if
$t_j$ is a solution so is $-t_j$.
The equation can be rewritten to determine $\lambda_1$ directly:
\begin{equation}
 \tan\left(\frac{\sqrt{3\lambda_1}}{2}\right)
= \frac{1}{\beta} 
\frac{3\sqrt{\lambda_1}} 
{\left(\lambda_1 -\frac{2}{\beta^2}\right)}
\label{eq:TriLambdaBeta}
\end{equation}

The account above suffices for calculating $\lambda_1$: see the further development following this paragraph.
The set of all eigenvalues and eigenfunctions are determined in \cite{McC04,FK14}.
The approach in McCartin~\cite{McC04} leads to a coupled pair of transcendental equations which we now
show is equivalent to our equation~(\ref{eq:transct}) as follows.
First our interest is in McCartin's `symmetric modes' and the eigenvalue is given by his~(5.3)
$$ \lambda_1=\frac{4}{27} \left(\frac{\pi}{r}\right)^2 (\mu^2 + \mu\nu +\nu^2) ,$$
where $r$ is the inradius of the triangle and $\mu$ and $\nu$ are described below.
They are found in terms of quantities $L$, $M$ and $N$ solutions of coupled  transcendental equations
given in McCartin's equations~(5.7) where, for our fundamental mode we can set $m=n=0$ and $M=N$ 
arriving at
\begin{eqnarray*}
2 (L-M) \tan(L) &=& \frac{3 r}{\beta} , \\
-(L-M) \tan(M) &=& \frac{3 r}{\beta} . 
\end{eqnarray*}
Since
$$ \tan(L-M)= \frac{\tan(L)-\tan(M)}{1+\tan(L) \tan(M)} $$
the coupled equations can be wtitten as a single equation for $(L-M)$
in which $\tan(L-M)$ is equal to a simple rational function of $(L-M)$.
Then setting $(L-M)=t\sqrt{3}/\beta$ in the single equation and using $r=1/\sqrt{3}$
we recover the transcendental equation given in  our equation~(\ref{eq:transct}).
McCartin's formula~(5.8) gives, for our fundamental mode,
$$\mu=\nu= -\frac{L-M}{\pi} $$
and hence
$$\lambda_1
=\frac{4}{9} \left(\frac{\pi}{r}\right)^2 \mu^2 
= \frac{4}{9} \left(\frac{1}{r}\right)^2 (L-M)^2 
= \frac{4}{3}(L-M)^2
$$
on using $r=1/\sqrt{3}$.
Finally, substituting $(L-M)=t\sqrt{3}/\beta$ gives
$$\lambda_1 = \frac{4 t^2}{\beta^2}$$
as previously given in equation~(\ref{eq:lambda1TriMK}).

There are several approaches to equation~(\ref{eq:TriLambdaBeta}) which yield information
(some of which may be in Lam{\'e}'s paper).
\begin{itemize}
\item It is easy to solve for $\beta$ as a function of $\lambda_1$ though considerations
of branches of solutions do enter:
\begin{equation}
    \beta(\lambda) = \left\{\begin{array}{ll}
        \frac{3+\sqrt{9 +8 \tan(\sqrt{3\lambda}/2)^2}}{2 \sqrt{\lambda} \tan(\sqrt{3\lambda}/2)}
        , & \text{for } 0\leq \sqrt{\lambda} < \frac{\pi}{\sqrt{3}}\\
        \frac{\sqrt{6}}{\pi}  
       ,  & \text{for } \sqrt{\lambda}=\frac{\pi}{\sqrt{3}} \\
        \frac{3-\sqrt{9 +8 \tan(\sqrt{3\lambda}/2)^2}}{2 \sqrt{\lambda} \tan(\sqrt{3\lambda}/2)}
       , & \text{for } \frac{ \pi}{\sqrt{3}} < \sqrt{\lambda}
        \end{array}\right.
\label{eq:betaLamTri}
\end{equation}
 Asymptotics are easily found. For example
 $$\beta(\lambda)\sim \frac{2\sqrt{3}}{\lambda} \qquad
 {\rm as\ }\ \lambda\rightarrow 0 .$$
 This is consistent with the general result that
 $$\lambda_1\sim \frac{|\partial\Omega|}{\beta|\Omega|}  \qquad
 {\rm as\ }\ \beta\rightarrow\infty .$$
 \item By differentiating equation~(\ref{eq:TriLambdaBeta})  one finds;
 \begin{equation}
 \frac{d\lambda}{d\beta}
 =  \frac{-12 \lambda (\beta^2 \lambda +2)}
 {\sqrt{3}\beta^4 \lambda^2 + \beta^2 (5\sqrt{3}+ 6\beta)\lambda + 4\sqrt{3}+12\beta} \  .
 \label{eq:TridLambda1Deta}
 \end{equation}
 Solving the initial-value problem starting from 
 $\lambda_1(\beta=0)= (4\pi^2)/3$ finds the $\lambda_1$ at larger values of $\beta$.
 For any positive starting value for $\lambda(0)$, $\lambda(\beta)$ decreases as
 $\beta>0$ increases.
 \end{itemize}
Asymptotics for $\beta$ large and for $\beta$ small can be found.
We have
$$ \lambda_1 \sim \frac{2\sqrt{3}}{\beta}\qquad {\rm for}\  \beta\rightarrow\infty. 
$$
Also
$$ \lambda_1 \sim\frac{4\pi^2}{3} - \frac{8\pi^2}{\sqrt{3}} \beta  \qquad {\rm for}\  \beta\rightarrow{0}. 
$$

\begin{figure}[hb]
\centerline{
$$
\includegraphics[height=5cm,width=6.5cm]{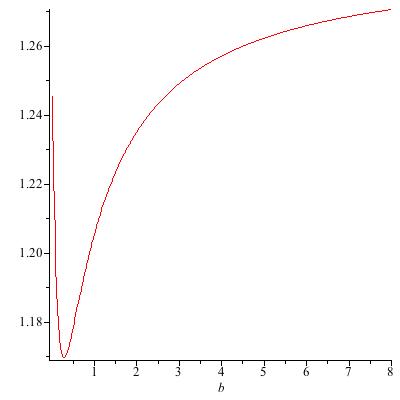}
$$
}
\caption{For domains with area $\sqrt{3}$, the ratio $\lambda_{1\Delta}/\lambda_{1\odot}$
plotted against $\beta$.
For domains both with the same area, but different from $\sqrt{3}$, the curves are different
but remain greater than 1}
\label{fig:equilatTriCirclambda1bNe0}
\end{figure}

\clearpage


\subsection{$Q_\steady$ for rectangle}\label{subsec:rectQs}

The pulsatile flow in a rectangular channel is treated in detail in this journal: see~\cite{WWS14}.
Our main concern is isoperimetric results appropriate to $Q_\steady$, the steady flow
(limit of the pulsatile flow as the frequency tends to zero) and, in  later subsection, with $\lambda_1$ which determines
the speed transients die out.

Consider the rectangle $(-a,a)\times(-b,b)$.
\smallskip

$Q_\steady$  is given by formula~(40) with (26) and (27) of~\cite{WWS14}.
There are different approaches, and earlier derivations.
See~\cite{ES65,Wa12b} and~\cite{MK94}. 
A very slight adaptation of material in~\cite{MK94} is as follows.

Let
\begin{equation}
 u_\steady
= \frac{a^2-x^2}{2} +\beta a -
\sum_{p=1}^\infty C_p \cosh\left( \frac{X_p y}{a}\right) \cos\left( \frac{X_p x}{a}\right)  ,
\label{eq:uSrect}
\end{equation}
with $X_p$ satisfying
$$ \cos(X_p)
= \frac{\beta}{a} X_p \sin(X_p)
$$
so that the boundary conditions at $x=\pm{a}$ are satisfied.
To satisfy the  boundary conditions at $y=\pm{b}$, we calculate that at $y=+b$
$$ u +\beta\frac{\partial u}{\partial y}
= \frac{a^2-x^2}{2} +\beta a -
\sum_{p=1}^\infty C_p\left(  \cosh\left( \frac{X_p b}{a}\right)+\frac{\beta}{a}\sinh\left( \frac{X_p b}{a}\right)\right)
 \cos\left( \frac{X_p x}{a}\right) .
 $$
 We find $C_p$ as follows. Write
 $$\frac{a^2-x^2}{2} +\beta a
 = \sum_{p=1}^\infty A_p  \cos\left( \frac{X_p x}{a}\right) \ \ 
 {\rm so} \ \
 A_p = 2\left(\frac{a}{X_p}\right)^3 \frac{\sin(X_p)}{(\beta\sin(X_p)^2 +a)} .
 $$
 Then
 $$ C_p
 = \frac{A_p}{ \cosh\left( \frac{X_p b}{a}\right)+\frac{\beta}{a}\sinh\left( \frac{X_p b}{a}\right)} .
 $$
 Integrating~(\ref{eq:uSrect}) gives
 $$
 Q_\steady
 = \frac{4}{3} a^3 b + 4\beta a^2 b -
 4 a^2 \sum_{p=1}^\infty C_p \frac{\sin(X_p)\sinh\left(\frac{b X_p}{a}\right)}{X_p^2} .
 $$

\begin{figure}[hb]
\centerline{
$$
\includegraphics[height=5cm,width=6.5cm]{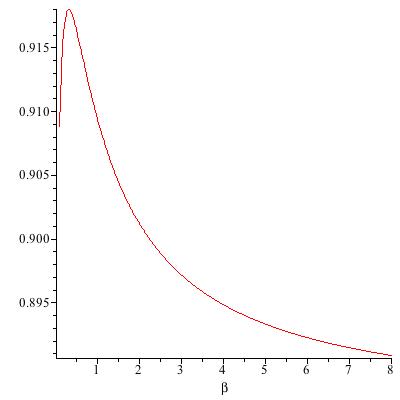}
$$
}
\caption{For domains with area $\pi$, the ratio $Q_\steady({\rm square})/Q_\steady({\rm circle})$
plotted against $\beta$.
For domains both with the same area, but different from $\pi$, the curves are different
but remain less than 1}
\label{fig:squareQratio}
\end{figure}

Figure~\ref{fig:squareQratio} indicates that,
for domains with the same area, $Q_\steady({\rm square})$ is
less than $Q_\steady({\rm circle})$.
The ratio shown there varies only by about 10\% over the whole range of $\beta$.
 $Q_\steady({\rm circle})$ is linear in $\beta$, and plots of $Q_\steady({\rm square})$
 are approximately linear in $\beta$.
Plots, for other rectangles also indicate that $Q_\steady$ is approximately linear in $\beta$.
This is in agreement with the statement in the Conclusion of~\cite{WWS14}.

\subsection{$\lambda_1$ for rectangles generally}\label{subsec:Rectangular} 

The function $u=\cos(\mu_X x)\cos(\mu_Y y)$ satisfies
$$ \Delta u + \lambda u = 0 ,\qquad {\rm with}\ 
\lambda = \mu_X^2 + \mu_Y^2 . $$
The Robin boundary conditions are satisfied if
\begin{equation}
\mu_X\tan(\mu_X a)=\frac{1}{\beta} , \qquad
 \mu_Y\tan(\mu_Y b)=\frac{1}{\beta} .
 \label{eq:muXmuY}
 \end{equation}
%
 A `dimensionless'  form of the equation is
  \begin{equation}
{\hat\mu} \tan({\hat c}{\hat\mu})=1,
 \qquad{\mbox{\rm where\ \ }} {\hat\mu}=\beta\mu,\ {\hat c}=\frac{c}{\beta} \  .
 \label{eq:rectmuchat}
 \end{equation}
 The transcendental equations have been widely studied, e.g.~\cite{BS73,MRS03,LWH15}.
 Numerical values, often used for checks, are given in Table 4.20 of 
 Abramowitz and Stegun, {\it Handbook of Mathematical Functions} (1964, Dover: 1965).
 
 We have an interest in the smallest positive solutions,
 $$0<\mu_X<\pi/(2 a),\qquad  0<\mu_Y<\pi/(2 b),\qquad  0<\mu<\pi/(2 c) .
 $$
 Because it is the square of $\mu$ which occurs in $\lambda_1$
 (as in equation~(\ref{eq:lambda1Rect})), we also define
 $\mu_{(2)}=\mu^2$.
 \begin{itemize}
 \item
 At fixed $c>0$, $\mu$ decreases as $\beta$ increases.
 \medskip
 
 For $\beta$ small, $\mu\sim\pi(1-\beta/c)/(2c)$\\
 for $\beta$ large,
 \begin{equation}
 \mu\sim{1/\sqrt{c\beta}} \qquad{\rm for}\ \beta\rightarrow\infty . 
 \label{eq:mucbetaInfty}
 \end{equation}
 
 \item
 At fixed $\beta>0$, $\mu$ decreases as $c$ increases:
 \begin{equation}
 \frac{d  \mu}{d  c}
 = - \ \frac{\mu ( 1+ \beta^2\mu^2)}{\beta +c (1+\beta^2\mu^2)} ,\qquad
 \frac{d  \mu_{(2)}}{d  c}
 = - \ \frac{2\mu_{(2)} ( 1+ \beta^2\mu_{(2)})}{\beta +c (1+\beta^2\mu_{(2)})} .
\label{eq:dmu2dc}
\end{equation}
The $\mu_{(2)}$ de can be written
$$ \frac{d  {\hat\mu}_{(2)}}{d {\hat c}} 
= - \ \frac{2{\hat\mu}_{(2)} ( 1+ {\hat\mu}_{(2)})}{1 +{\hat c} (1+{\hat\mu}_{(2)})} .
$$
Also $\mu(c)$ is convex in $c$:
$$\frac{d ^2 \mu}{d  c^2}
 =  \ \frac{2\mu ( 1+ \beta^2\mu^2)(2\beta^3 \mu^2 +c \beta^4 \mu^4 +2 c \beta^2 \mu^2 + \beta + c)}
 {(\beta +c (1+\beta^2\mu^2))^3} .
 $$
Continuing calculations we see that $\log(\mu(c))$ is convex in $c$:
 \begin{eqnarray*}
 \frac{d^2 \log(\mu(c))}{d c^2}
 &=& \frac{ \mu(c)\mu''(c)-\mu'(c)^2}{\mu(c)^2}\\
 &=&  \ \frac{(1+ \beta^2\mu^2)(3\beta^3 \mu^2 +c \beta^4 \mu^4 +2 c \beta^2 \mu^2 + \beta + c)}
 {(\beta +c (1+\beta^2\mu^2))^3} .
\end{eqnarray*}
On the basis of having calculated a few more higher derivatives  $\frac{d ^j \mu}{d  c^j}$  of $\mu(c)$ with respect to $c$
there are indications that $\mu(c)$ may be completely monotone.
(Any completely monotone function is both nonincreasing and logconvex.)
In Lemma 1 of \S\ref{subsec:Proofexplicit} we show that $c(\mu)$ is completely monotone.
\medskip
 
For $c$ small, $\mu\sim{\pi/(2\sqrt{c\beta})}$.\\
For $c$ large, $\mu\sim{\pi/(2({c+\beta}))}$.
\end{itemize}
 
The fundamental Robin eigenvalue for the rectangle is,
with $\mu_X$ and $\mu_Y$ each the smallest root of their defining transcendental equations,
\begin{equation}
\lambda_1 = \mu_X^2 + \mu_Y^2 .
\label{eq:lambda1Rect}
\end{equation}
 
We have, for a fixed rectangle, the following asymptotics,
where we write $a=r h$ and $b=h/r$:
\begin{eqnarray}
\lambda_1
&\sim& \frac{\pi^2}{4}\left( 
\frac{1}{a^2}+ \frac{1}{b^2}- \frac{2\beta}{a^3}- \frac{2\beta}{b^3} \right)\qquad
{\rm as\ \ } \beta\rightarrow{0}
\label{eq:lambda1Rectb0}\\
&\sim& \frac{\pi^2}{4h^2}\left( 
\left( \frac{1}{r\left(1+\frac{\beta}{h r}\right)}\right)^2 +
\left( \frac{r}{\left(1+\frac{\beta r}{h}\right)}\right)^2 
\right) \qquad
{\rm as\ \ } \beta\rightarrow{0}
\label{eq:lambda1Rectb0r}\\
 \lambda_1
&\sim& \frac{|\partial\Omega|}{\beta|\Omega|} 
=\frac{1}{\beta}\left( \frac{1}{a}+ \frac{1}{b}\right)
=\frac{1}{\beta h}\left( \frac{1}{r}+ r\right) \qquad
{\rm as\ \ } \beta\rightarrow\infty .
\label{eq:lambda1RectbInf}
\end{eqnarray}

In the case of a square $\mu_X=\mu_Y=\mu$ and 
$\lambda_{1\square}=\lambda_1({\rm square})=2\mu^2$ so the
transcendental equation for $\mu$ is readily re-written in terms of $\lambda_1$.
Also, with $a=b=h$,
$$\frac{d \lambda_{1\square}}{d \beta}= \frac{-2\lambda_{1\square}}{\beta+h(1+\beta^2 \lambda_{1\square}/2)} ,
\qquad
\lambda_{1\square}(\beta=0)= \frac{\pi^2}{2 h^2} ,
$$
which, for a square, corresponds to an equation~(\ref{eq:TridLambda1Deta})
for an equilateral triangle. 
In both cases one sees that $\lambda_1(\beta)$ decreases as $\beta$ increases.
For a square one can show that $\lambda_{1\square}(\beta)$ is log-convex in $\beta$.
The  terms in 3rd derivatives and higher involve both signs.
For more on $\lambda_{1\square}$ see the next subsection~\S\ref{subsec:Square}.
For rectangles more generally there does not appear to be a single transcendental equation,
or differential equation, for $\lambda_1(\beta)$ but coupled equations involving $\mu_X$ and $\mu_Y$
appear to be needed.

\begin{figure}[h]
\centerline{
$$
\includegraphics[height=5cm,width=6.5cm]{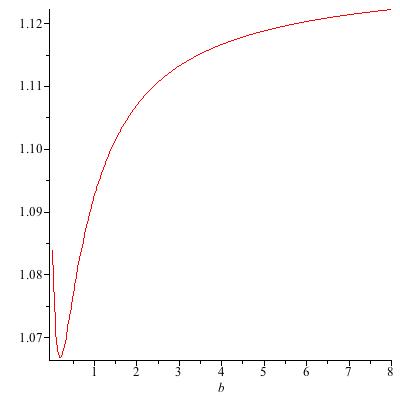}
$$
}
\caption{For domains with area $\pi$, the ratio $\lambda_1({\rm square})/\lambda_1({\rm circle})$
plotted against $\beta$.
For domains both with the same area, but different from $\pi$, the curves are different
but remain greater than 1}
\label{fig:lambda1SquareCircplot}
\end{figure}

It is generally true that amongst domains with the same area, that which has the smallest fundamental eigenvalue is a disk.
A numerical illustration of this fact, for the square, is given in Figure~\ref{fig:lambda1SquareCircplot} .

Various checks of the results shown in Figure~\ref{fig:lambda1SquareCircplot} are possible.
When $\beta=0$, 
$\lambda_1({\rm circle})= j_0^2 \approx 5,783$ with $j_0$ the smallest zero of $J_0$,
and $\lambda_1({\rm square})= 2\pi$, so the ratio is approximately 1.086.

\medskip
 
We also compare $\lambda_1$ for rectangles of the same area with each other.
 See \S\ref{subsubsec:PolygonalbGt0},
in particular,  Figure~\ref{fig:lambda1Rectb1SquareMin}.
In this connection,
further development from equation~(\ref{eq:lambda1Rect}) is given in Theorem~\ref{thm:thm3} in \S\ref{sec:rectlambda1}.

\subsection{Inequalities for $\lambda_{1\square}$ for a square, side $2h$}
\label{subsec:Square}

The formula for $\mu_\square^2$ can be rewritten
\begin{equation}
\frac{1}{\beta}
= \sqrt{\frac{\lambda_{1\square}}{2}}\tan\left(\frac{\sqrt{|\Omega|}}{2}\frac{\lambda_{1\square}}{2}\right) .
\label{eq:betaLamSq}
\end{equation}
(Compare this with equation~(\ref{eq:betaLamTri}).)

Several items are immediate consequences of items in the previous subsection.
Equation~(\ref{eq:mucbetaInfty}) gives
$$ \lambda_\square = 2\mu^2
\sim\frac{2}{h\beta} = \frac{|\partial\Omega|}{\beta|\Omega|} \qquad
{\rm as\ }\ \beta\rightarrow\infty .
$$

The following inequalities on $\tan(x)$
\begin{equation}
\frac{8 x}{\pi^2-4 x^2}
< \tan(x)
< \frac{x \pi^2}{\pi^2-4 x^2}
\label{eq:A1}
\end{equation}
are established in~\cite{ZLGQ}.

Upper bounds on the tan function lead to lower bounds on $\mu(r)$ while
lower bounds on the tan function lead to upper bounds on $\mu(r)$.
Inequalities~(\ref{eq:A1})
lead to
\begin{equation}\mu_{\rm LB}(c,\beta)=
{\frac {\pi }{\sqrt {c \left( 4\,c+{\pi }^{2}\beta \right) }}}
<\mu(c,\beta)
<\mu_{\rm UB}(c,\beta)= \frac{\pi}{4\beta}\, \left( -1+\sqrt {1+4\,{\frac {\beta}{c}}} \right)  .
\label{eq:LBUBmu}
\end{equation}
In particular, the formula for $\mu_{\rm LB}(h,\beta)$ gives
$$\lambda_\square > \frac{2}{\beta h+ 4 h^2/\pi^2} . $$

Other bounds can be established. For example, for all $r>0$,
$$\lambda_\square >  \frac{2}{h\beta} -\frac{2}{\beta^2} - \frac{4 r (r^2+1)}{\beta h(r^4+1)} ,\qquad
{\rm so\ }\ \lambda_\square >  \frac{2}{h\beta} -\frac{2}{\beta^2} . 
$$
This follows from equation~(\ref{eq:VxVymuSq})
on using that ${\cal V}_x>0$ and ${\cal V}_y>0$.

\subsection{The modulus of asymmetry for the equilateral triangle and for rectangles}

See~\cite{Fi14} for some calculations of the modulus of asymmetry 
(see~\S\ref{subsec:modulusAsym}) for various domains.
The interest here is in how it might relate to quantitative estimates of the difference between
the functionals for $\Omega$ and for the disk of the same area $\Omega^*$.
\clearpage

\section{Steady flows}\label{sec:Steady}

\subsection{The main theorem}\label{subsec:isoperimetricQSteady}

There have been many studies of these, sometimes oriented towards exact solutions,
sometimes numerical approximations, and sometimes checking theory against experiment:
see e.g.~\cite{BTT09,TH11,Wa12,Wa12b} and~\cite{Wa03}.

For given $\Omega$ and $\Delta{p}>0$ there is a 
unique solution for $u$ and it is positive in $\Omega$ (Theorem 2.3 of~ \cite{KM93},
but well known before this).
The solutions at different $\slipparameter$ vary monotonically with $\slipparameter$ 
(as shown in Theorem 2.4 of ~\cite{KM93}). 
The effects of varying the domains $\Omega$ are more difficult to study:
there is monotonicity with domain inclusion when $\slipparameter=0$ 
(Theorem 2.6 of~ \cite{KM93}) 
but not generally when $\slipparameter\ne{0}$.

We now reference the proof of the first sentence in the Abstract.
The area of $\Omega$ is denoted by $|\Omega|$.
When $\slipparameter=0$ the result is called the St Venant inequality.
The proof techniques for non-zero $\slipparameter$ differ from the `symmetrisation' proof techniques
used for the earlier proofs for $\slipparameter=0$.
The pure mathematical papers write that they
denote by $\Omega^*$ the `ball' with the same `volume' as $\Omega$: in our
2-dimensional case, $\Omega^*$ is the (circular) disk of radius $\sqrt{|\Omega|/\pi}$.

\medskip
\begin{theorem}
\label{thm:thm1}
 For every $\slipparameter > 0$, for every bounded Lipschitz set $\Omega\subset{R^N}$, the following inequality holds
\begin{equation}
Q_{\steady}(\Omega,\slipparameter) \le Q_\steady({\Omega^*},\slipparameter) . 
\label{eq:thm1}
\end{equation}
\end{theorem}
\medskip
\par\noindent
This is Theorem 1.1 of~\cite{BG15}.
(Caution. The $\beta$ of~\cite{BG15} is the reciprocal of our $\slipparameter$.)
The expression on the right hand side of equation~(\ref{eq:thm1})
is given by equation~(\ref{eq:circleQinf}) with
$a=\sqrt{|\Omega|/\pi}$.

Formulae for $Q_\steady(\beta)$ for various shapes are given in~\cite{MK94}.
Other inequalities, for general $\Omega$, are given in~\cite{KM93}.
Related work is presented in~\cite{BB13}.

\subsection{$Q_\steady$: circle compared with equilateral triangle and square}
\label{subsec:QscircleTriSq}

The formula for $Q_\steady({\rm circle})$ is given in equation~(\ref{eq:circleQinf}).
Figures~\ref{fig:equilatTriQratio} (for the triangle)
and \ref{fig:squareQratio} (for the square) provide examples consistent with Theorem~\ref{thm:thm1}.
\clearpage

\section{Transient Flows}\label{sec:Transient} 

\subsection{An isoperimetric result} \label{subsec:isoperimetricLambda1}

The second sentence of the Abstract follows from the representation~(\ref{eq:uStartingSeries}),
noting that $\lambda_1$ is the least of the eigenvalues, and the following theorem.

\begin{theorem}
\label{thm:thm2}
For every $\beta > 0$, for every bounded open Lipschitz set $\Omega\subset{R^N}$
\begin{equation}
 \lambda_1(\Omega,\slipparameter) \ge \lambda_1({\Omega^*},\slipparameter) .  
\label{eq:thm2}
\end{equation}
with equality if and only if $\Omega=\Omega^*$.
\end{theorem}
\medskip
\par\noindent
This is the $N=2$, $q=2$ case of Theorem 2.1 
of~\cite{BG15}.
(The paper~\cite{BG15} titles Theorem 2.1 as
{\rm A family of Faber-Krahn inequalities}.) 
The theorem was first proved, in the $N=2$ case in~\cite{Bo86CR} (see also~\cite{Bo86CV,Bo88},
and the general $N$ case in~\cite{Da06}.
Related work is presented in~\cite{BD10,BG10,BG15eig,Da00,Da06}.
The expression on the right hand side of equation~(\ref{eq:thm2})
is given by equation~(\ref{eq:lambda1gamma}) with
$a=\sqrt{|\Omega|/\pi}$.

\subsection{$\lambda_1$: circle compared with equilateral triangle and square}
\label{subsec:lambda1circleTriSq}

Figures~\ref{fig:equilatTriCirclambda1bNe0}(for the triangle)
and \ref{fig:lambda1SquareCircplot} (for the square) provide examples consistent with Theorem~\ref{thm:thm2}.

\clearpage

\section{Illustrative examples of, and notes on, time-dependent flows}\label{sec:IllustrativeTimeDep}


\subsection{$\beta=0$ and starting flows}
\label{subsec:QtcircleTSq}

\begin{figure}[hb]
\centerline{
$$
\includegraphics[height=4cm,width=6.5cm]{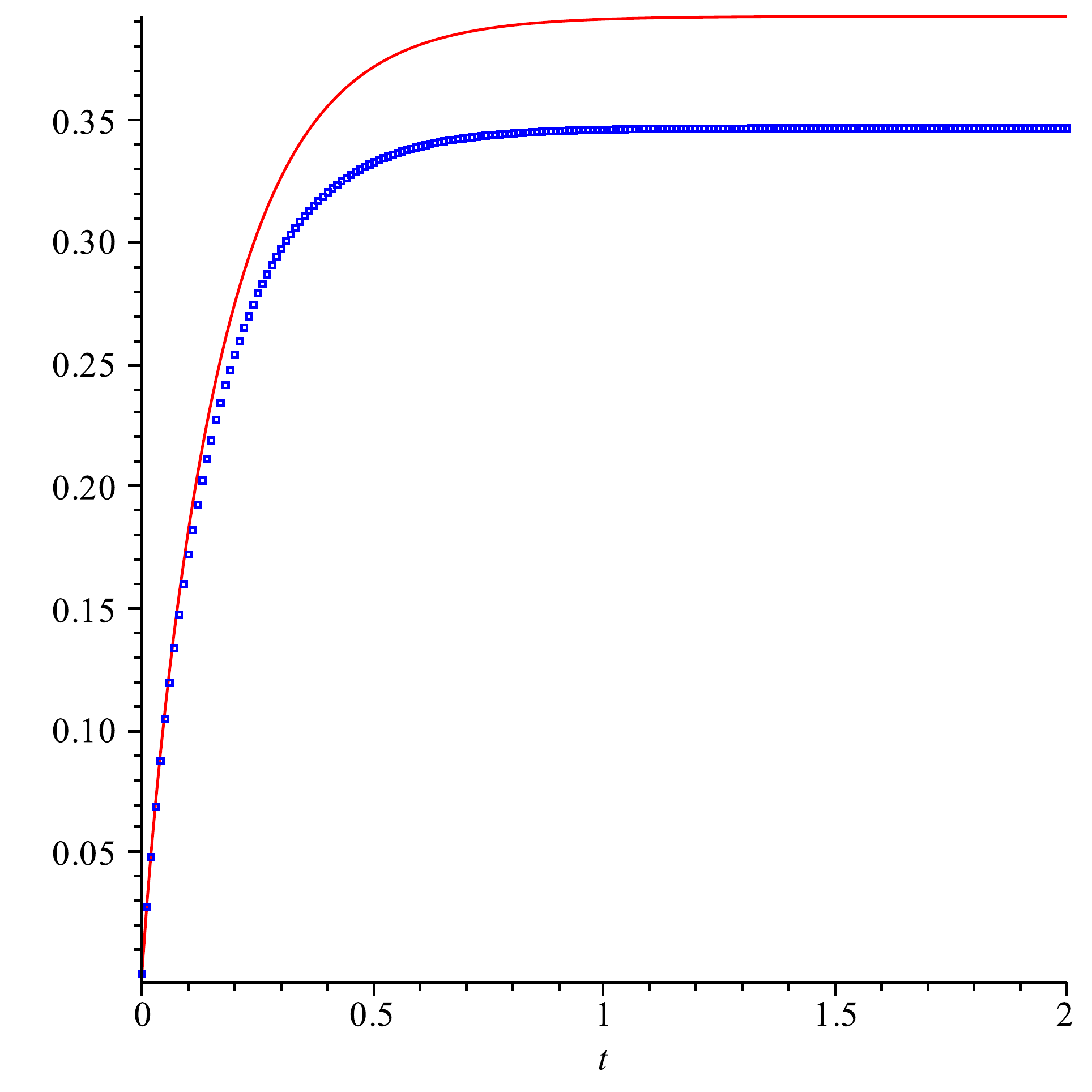}
\includegraphics[height=4cm,width=6.5cm]{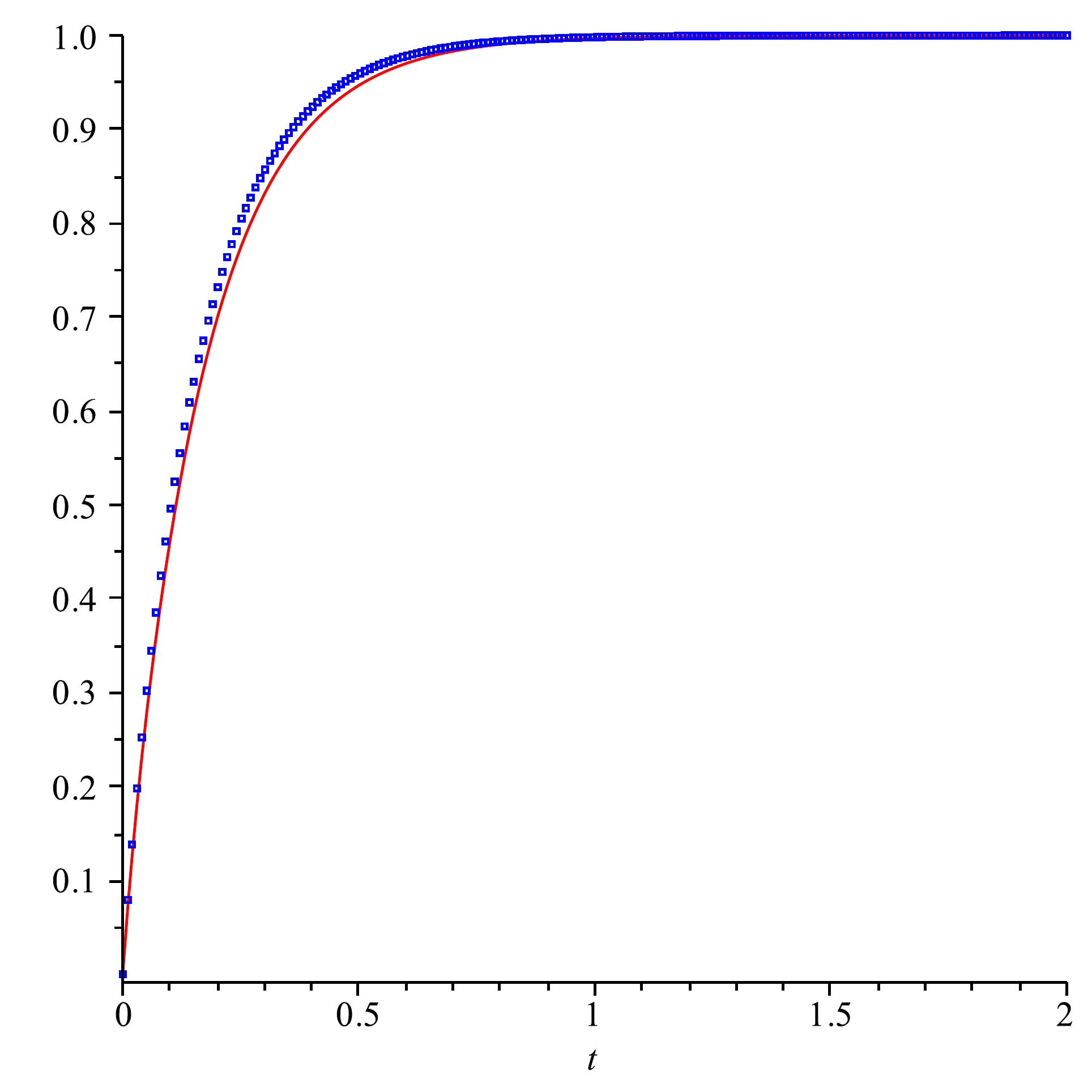}
$$
}
\caption{The left hand plots are of $Q(t)$:
the right hand plots are of $Q(t)/Q(\infty)$.
The plots for the circular cross section are shown by the solid (red) curve:
those for the square cross section are shown dotted (blue) with a small square symbol.}
\label{fig:zerocs}
\end{figure}

The results of Theorems 1 and 2 have been known in the $\beta=0$ case
since 1949 and 1923 respectively.
We suspect that graphical illustration of the kind shown in Figure~\ref{fig:zerocs}
will be available in older literature, perhaps a century ago.
The separation-of-variables, eigenfunction expansions are straightforward to find
in several domain shapes, circles and rectangles being common examples.
We have already cited~\cite{Ba67} for the circle while
the rectangle is treated, for example, in~\cite{Er03},
though we expect the first publication of these might well be over a century ago.
In Figure~\ref{fig:zerocs}
we show plots of $Q(t)$ for a circle of unit radius and a square with sides $\sqrt{\pi}$.

The left hand plots illustrate that the disk has the larger $Q(\infty)$.
The right hand plots illustrate that the approach to $Q(\infty)$ is slower
for the circular cross section than for others.

It is also easy to plot $Q(t,{\rm square})/Q(t,{\rm circle})$.
The values in~\cite{PoS51} can be used to check the behaviour at very large time as,
for domains with area $\pi$,
$Q_\steady({\rm circle})= \pi/8\approx{0.3927}$ and
$Q_\steady({\rm square})\approx{0.3469} $.

\subsection{$\beta>0$: comparing $\lambda_1$}

The results are less striking for large $\beta$.
In particular, when $\beta$ is very large the situation is essentially like that when $u_n$ is an appropriate constant around the boundary,
$u$ is constant across $\Omega$ and
$Q(t)$ is approximately $t|\Omega|\Delta{p}$,
and there is no dependence on the shape of $\Omega$
(at the lowest order in the approximation for $\beta$ large).

\subsection{A few more results valid for $\beta>0$}
\label{subsec:otherBetaGt0}

\subsubsection{Other checks on $\phi_n$}

See~\cite{PoS51} p107:
$$|\Omega| 
= \sum_{n=1}^\infty
\frac{\left( \int_\Omega \phi_n \right)^2}{\int_\Omega \phi_n^2} .
$$
Similarly the polar moment of inertia about the centroid satisfies
$$ I_c
= \sum_{n=1}^\infty
\frac{\int_\Omega |z-z_c|^2\phi_n(z) \int_\Omega \phi_n }{\int_\Omega \phi_n^2} .
$$

\subsubsection{$Q_\steady$}

$$
\leqno{\bf\rm\bf Conjecture} \qquad\qquad\qquad
\frac{Q_\steady(\Omega^*,\beta)-Q_\steady(\Omega,\beta)}
{Q_\steady(\Omega^*,\beta)} \rightarrow 0\qquad
{\rm as} \ \ \beta\rightarrow\infty .
$$

\subsubsection{The decay rate $\lambda_1$ when $\beta>0$}

Some monotonicity results for the principal eigenvalue of the generalized Robin problem
are established in~\cite{GS05}.
Their $\alpha=-1/\beta$; their $\lambda$ has the opposite sign to ours,
and it is our sign we use in this report.
Before treating monotonicity aspects we repeat items they report:
The eigenvalue $\lambda_1(\beta)$ is simple, and that 
an eigenfunction $\phi_1(.,\beta)$ can be chosen with a single sign 
and normalized by letting the integral of its square over $\Omega$ being 1.
In addition, as $\beta\rightarrow{0}$,
$\lambda_1(\beta)$ converges to the principal eigenvalue for the Dirichlet problem.
Also
$$ \lim_{\beta\rightarrow\infty}{\lambda_1(\beta)}{\beta}
= \frac{|\partial\Omega|}{|\Omega|} .
$$
There are partial monotonicity results.
\begin{itemize}
\item
Let $B\subset{R^n}$ $n\ge{2}$ be a ball and $\Omega\subset{B}$.
If $\alpha<0$ then 
$\lambda_1(\Omega_1,\alpha)>\lambda_1(B,\alpha)>0$.
\item
Again, let $\alpha<0$. 
If $\Omega\subset{R^n}$, $n\ge{2}$, is a convex domain that contains a ball $B$, then
$$\lambda_1(B,\alpha){\ge}\lambda_1(\Omega,\alpha)>0.
$$ (See Theorem 2 and Corollary 3 of ~\cite{GS05}.)
\item
Let $\alpha<0$ and $\Omega_1\subset{B}\subset{\Omega_2}$  
where $B\subset{R^n}$ $n\ge{2}$ is a ball and
$\Omega_1$, $\Omega_2$ are convex domains, we have that 
$\lambda_1(\Omega_1,\alpha){\ge}\lambda_1(\Omega_2,\alpha)>0$ 
\end{itemize}
\medskip

\subsection{Further results applicable when $\beta=0$}
\label{subsec:furtherBeta0}

\subsubsection{Isoperimetric results for some geometric functionals}

Strictly speaking this subsubsection doesn't require $\beta=0$.
However, at present our applications of the geometric functionals in connection
with $Q_{\rm steady}$ and $\lambda_1$ are just for the case $\beta=0$.

In this version of these notes,
the other geometric domain functionals we consider are, the perimeter $|\partial\Omega|$,
the polar moment of inertia $I_c$
and $B$ defined -- for star-shaped domains --  by equation~(\ref{eq:Bdef}).
\begin{equation}
 B(z_0)= \int_{\partial\Omega} \frac{ds}{n\cdot(z-z_0)} \qquad
{\rm and\ }\qquad
B= \min_{z_0\in\Omega} B(z_0) .
\label{eq:Bdef}
\end{equation}

For domains in the plane, see~\cite{PoS51} for Polya's proof that 
\begin{equation}
\frac{2I_c}{\pi} \leq \bigl(\frac{|\partial\Omega|}{2\pi}\bigr)^4
\label{eq:ioL4} \end{equation}
Equality is attained only for disks.
(Hadwiger 
 gave a different proof, one which generalizes to $R^n$.
For references, see~\cite{Ke06}.)
We remark
that a combination of this~(\ref{eq:ioL4}) and another
easier-to-prove inequality, represents a refinement of the
classical $(|\Omega|,|\partial\Omega|)$, area-perimeter, isoperimetric inequality:
\begin{equation}
\bigl(\frac{|\Omega|}{\pi}\bigr)^2 \leq \frac{2I_c}{\pi}
\leq\bigl(\frac{|\partial\Omega|}{2\pi}\bigr)^4 .
\label{ineq:AIcL}
\end{equation}

We first note that connectedness is essential for
inequality~(\ref{eq:ioL4}).
For, if we let the domain 
be the union of two
equal disjoint disks, symmetically placed either side of the
origin so that the centroid is at the origin, we have a counterexample.
By taking the components further apart we can
increase $I_c$ indefinitely, while $|\partial\Omega|$ stays fixed.
However, on
joining the disks by a straight line the perimeter increases and
this dumbell shaped domain does satisfy
inequality~(\ref{eq:ioL4}).

Polya's proof begins with conformal mapping, and requires the domain 
to be simply connected.
Inequality~(\ref{eq:ioL4})  can be shown to be true for multiply-connected
ones too: see~\cite{Ke06}.

\subsubsection{Physical functionals: isoperimetric results applicable when $\beta=0$}

It was conjectured by Polya and Szego that among sets with given torsional rigidity
$Q_{\rm steady}$,  balls minimize $\lambda_1$. specifically
\begin{equation}
Q_{\rm steady} \lambda_1^2 \ge \frac{\pi}{8} j^4 . 
\label{eq:KJ}
\end{equation}
This is proved in~\cite{KJ78,KJ75}.

The Payne-Rayner inequality, proved in~\cite{PR72},
s relevant to estimating $Q_1(\beta=0)$.
It states
$$ Q_1 \ge \frac{4\pi}{\lambda_1^2} . $$
Furthermore, there is equality when $\Omega$ is a disk.
We remark that for the unit disk
$$ Q_n(\rm disk) = \frac{4\pi}{j_n^4} , $$
where $j_n$ is the $n$-th positive zero of $J_0$.

There are many more isoperimetric inequalities.
With $B$ defined -- for star-shaped domains --  by equation~(\ref{eq:Bdef}),
\cite{PoS51} p93 and p94 give
$$ Q_{\rm steady} \ge \frac{|\Omega|^2}{4 B}, \qquad
\lambda_1 \le \frac{ j^2 B}{2|\Omega|} .
$$
The first of these is an equality when $\Omega$ is an ellipse:
see~\cite{PoS51} p262.
Other inequalities involving $B$ include
$$\lambda_1 \le j^2 \frac{B}{2|\Omega|} , $$
(see~\cite{PoS51}) which gives an equality when $\Omega$ is a disk.
The functional $B$ occurs in an inequality involving $Q_1$,
$$ Q_1 \le \frac{4 B}{\lambda_1^2} . $$
See~\cite{CS78,He88}.

\subsubsection{Some other results when $\beta=0$ and $\Omega$ is convex}

When $\Omega$ is convex, 
\begin{itemize}

\item the steady solution, torsion function,  has a square root which is concave;

\item the principal eigenfunction $\phi_1>0$ is such that $\log(\phi_1)$ is concave.

\end{itemize}
Let $\Omega_0$ and $\Omega_1$ be convex, and define the Minkowski sum
$$ \Omega_t= (1-t)\Omega_0 + t \Omega_1 . $$
\begin{itemize}

\item $Q_{\rm steady}(0,\Omega_t)^{1/4}$ is concave.
See~\cite{Bo85}.

\item $\lambda_1(0,\Omega_t)^{-1/2}$ is  concave.

\end{itemize}

\section{Symmetrisation}\label{sec:Symmetrisation}

\noindent{\bf Theorem}~\cite{PoS51}
{\it Symmetrization of a plane domain with respect to a line leaves the
area unchanged, but decreases the perimeter, polar moment of inertia about the centroid.}\\
See~\cite{PoS51}p6. In connection with moments of inertia, Steiner symmetrization
about the line through the centroid reduces the moment associated with the
direction perpendicular to the line, and leaves unchanged that 
associated with the direction of the line.

See~\cite{PoS51} for the use of symmetrisation in proving
\begin{itemize}
\item for $\beta=0$, $Q_\steady$ increases,
\item for $\beta=0$, $\lambda_1$ decreases.
\end{itemize}
The methods are only applicable in the case $\beta=0$.
See also the Appendix to this paper.

\section{Overview: Perturbations from circular domains}\label{sec:Perturbations}

Theorems~\ref{thm:thm1} and~\ref{thm:thm2} say nothing 
about by how much the quantities for some cross-section  differ from 
those of a circular channel of the same cross-sectional area.
Furthermore for domains $\Omega$ with given area $|\Omega|$,
Theorem~\ref{thm:thm1} is an upper bound on $Q_{\rm steady}(\Omega)$,
while Theorem~\ref{thm:thm2} is an lower bound on $\lambda_1(\Omega)$.
Also of interest are bounds on the other side with various forms of
geometric restrictions.
And there are many other improvements that are possible,
and many more that have been conjectured but remain unproven.
For example, when $\beta=0$, \cite{PoS51}p112 report a result
(of Nicolai) from the 1920s that
\begin{equation}
Q_{\rm steady}
\le \frac{ I_{\rm max} I_{\rm min}}{I_{\rm max}+I_{\rm min}} ,
\label{eq:Nicolai}
\end{equation}
where $I_{\rm max}$ and $I_{\rm min}$ are the principal moments of
inertia about the centroid of $\Omega$.
Inequality~(\ref{eq:Nicolai}) is an equality for any ellipse.
With $I_c=I_{\rm max}+I_{\rm min}$ from this we readily see
$$Q_{\rm steady}
\le \frac{ I_{\rm max} I_{\rm min}}{I_{\rm max}+I_{\rm min}} 
\le \frac{I_c}{4} ,
$$
with the latter inequality an equality only for disks.

This section, and the next two, are primarily intended as a literature survey 
of items which seem relevant to 
these topics.
There is an historical component to this. 
Rayleigh supported his conjecture on what has become known as the
Faber-Krahn inequality with Fourier series approximations appropriate to nearly circular domains: 
this is described in~\cite{PoS51}.
Such asymptotics are the topic of~\S\ref{subsec:nearCirc}.
Before this we treat,
in~\S\ref{sec:Measuring},  methods less dependent on the perturbations being
small.
In~\S\ref{subsec:modulusAsym} we indicate an area-based measure
of asymmetry, one of many
used in improving isoperimetric inequalities.
In \S\ref{subsec:deficit} we review bounds in terms of the 
classical geometric isoperimetric deficit.
After these we consider nearly circular domains,
first for the case $\slipparameter=0$ in \S\ref{subsec:nearCircBeta0},
then, in work that hasn't previously been published,
 for ellipses with $\slipparameter>0$ in \S\ref{subsec:nearCircBetaNe0},

\section{Measuring the effect of departures from circular}\label{sec:Measuring}

\subsection{A measure of asymmetry}
\label{subsec:modulusAsym}

Several decades ago Fraenkel defined a modulus of asymmetry
$$ \alpha(\Omega) = {\rm min}_{z\in\Omega} \frac{|\Omega\setminus\Omega^*(z)|}{|\Omega|} ,
$$
where $\Omega^*(z)$ is the disk centred at $z$ with the same area as $\Omega$.
See~\cite{Fr08,Fi14}. Of course $0\le\alpha(\Omega)\le{1}$.
Since $\Omega$ and $\Omega^*(z)$ have the same area
$$|\Omega\setminus\Omega^*(z)|=|\Omega^*(z)\setminus\Omega | .$$
Note that some authors, perhaps most, use a symmetric difference, and hence have a modulus which is
twice that of the definition above.

\medskip\par\noindent{\bf $\alpha(\Omega)$ for ellipses $(x/a)^2 + (y a)^2=1$, $a\ge{1}$.}

The asymmetry $\alpha(\Omega)$ for an ellipse has been calculated several times,
e.g.~\cite{Fi14}  but this might be susceptible to improvement, absolute values used as needed, as
I would expect  $\alpha({\rm ellipse}(a))=\alpha({\rm ellipse}(1/a))$,
and we would expect an even function of $a$:
\begin{align*}
\alpha({\rm ellipse})
&=\frac{2}{\pi}\left( {\rm arcsin}\left(\frac{a}{\sqrt{1+a^2}}\right) -
 {\rm arcsin}\left(\frac{1}{\sqrt{1+a^2}}\right) \ \right) , \\
 &=\frac{2}{\pi}\left( {\rm arctan}(a)-
 {\rm arctan}(\frac{1}{a}) \ \right)
=\frac{2}{\pi} {\rm arctan}\left( \frac{a-1/a}{2}\right) ,\\
 &\sim
\frac{1}{\pi}\left(2 (a-1) - (a-1)^2 +\frac{1}{3} (a-1)^3 + \dots \ \right)\qquad{\rm as\ } a\rightarrow{1} .
\end{align*}
These expressions can be written in terms of $\epsilon$:
$$\epsilon
=\frac{a^2-a^{-2}}{a^2+a^{-2}} 
\sim 2(a-1) \qquad {\rm\ as\  }\ a\rightarrow{1} ,
$$
As
$$\left( a- \frac{1}{a}\right)^2 = a^2 +a^{-2} -2
=\frac{2}{\sqrt{1-\epsilon^2}} -2 ,
$$
we have
$$ \alpha({\rm ellipse})
=\frac{2}{\pi} {\rm arctan}\left(\sqrt{
\frac{1}{2}\left( \frac{1}{\sqrt{1-\epsilon^2}} -1\right)}
\right) .
$$
The ellipse is also treated in~\cite{HHW}, p89.

\subsubsection{Improving the isoperimetric inequality for $I_c$ in terms of $\alpha$}

This subsubsection is here in spite of the fact that,
with the inequalities on $Q_{\rm steady}$ and $\lambda_1$ presented in this note
to date, very few have any connection with $I_c$.
However, the section is included as it illustrates, perhaps with the easiest
domain functional for the purpose, that improvements to isoperimetric inequalities
are possible.

One domain functional for which the isoperimetric inequality can be easily improved using this modulus 
of asymmetry is $I_c$, the polar moment of inertia about the centroid.
The isoperimetric inequality is
$$  I_c(\Omega^*)=\frac{|\Omega|^2}{2\pi} \le I_c(\Omega)\ .$$
In 1986 (see~\cite{HHW}, p87) the first author proved
\begin{equation} 
  I_c(\Omega^*)\left( 1 + 2  \alpha(\Omega)^2\right ) \le I_c(\Omega) ,
\label{eq:Icasymm}
\end{equation}
which we have now seen, with a different proof, at page 198 of~\cite{HN92}.\\

\medskip

Consider polar moments about any point $z$ of $\Omega$, and denote the disk centred $z$ with the same area
as $\Omega$ by $\Omega^*_z$. Now, with the integrals being w.r.t area in the $z'$ variables
(i.e. $dx' dy'$)
$$ I_z(\Omega)= \int_\Omega |z-z'|^2 , \qquad
I_c(\Omega)= {\rm min}_{z\in\Omega} \ I_z(\Omega) \qquad
{\rm and}\ \   \int_{\Omega^*_z} |z-z'|^2 = \frac{|\Omega|^2}{2 \pi} .$$
Here $I_c$ denotes the moment about the centroid, $z_c$, and the Parallel-Axes Theorem 
(one of several ingredients in the proof in~~\cite{HN92}) states:
$$ I_z(\Omega) = I_c(\Omega) + |z- z_c|^2\, |\Omega| .
$$

Our proof of inequality~(\ref{eq:Icasymm}) is as follows.
Since $|\Omega|=|\Omega_z^*|$,
\begin{equation}
 |\Omega\setminus\Omega|=|{\Omega_z^*}\setminus\Omega| .
\label{eq:eqArea}
\end{equation}
We have
$$
\int_\Omega  |z-z'|^2 +  \int_{\Omega^*_z} |z-z'|^2
= \int_{\Omega\cup{\Omega^*_z}}  |z-z'|^2 + \int_{\Omega\cap{\Omega^*_z}}  |z-z'|^2 .$$
This is
 $$ I_z(\Omega) + \frac{|\Omega|^2}{2\pi}  
 = \int_{\Omega\cup{\Omega^*_z}}  |z-z'|^2 + \int_{\Omega\cap{\Omega^*_z}}  |z-z'|^2 .$$
Both the tems on the right can be bounded using the isoperimetric inequality, so
\begin{align*}
 I_z(\Omega) + \frac{|\Omega|^2}{2\pi}  
&\ge \frac{|{\Omega\cup{\Omega^*_z}} |^2 + |{\Omega\cap{\Omega^*_z}} |^2}{2\pi}\\
&= \frac{( |\Omega|+|{\Omega^*_z}\setminus\Omega|)^2+ (|\Omega|-|\Omega\setminus{\Omega^*_z}|)^2}{2\pi}\\
&=  \frac{ 2|\Omega|^2 + 2|\Omega\setminus{\Omega^*_z}|^2}{2\pi} .
\end{align*}
In the last line we have used equation~(\ref{eq:eqArea}).
We are free to choose $z$, so take $z=z_c$, the centroid.
Rearranging the last inequality gives
$$ I_c(\Omega)
\ge \frac{ |\Omega|^2}{2\pi}\left( 1 + 2 \frac{|\Omega\setminus{\Omega^*_{z_c}}|^2}{|\Omega|^2}\right)
\ge\frac{ |\Omega|^2}{2\pi}\left( 1 + 2\alpha(\Omega)^2\right) ,
$$
as $\alpha(\Omega)$ is found from minimizing over the centres $z$ of the disks.
This establishes the inequality~(\ref{eq:Icasymm}).
It is perhaps remarkable that repeated uses of the isoperimetric inequality yield an improvement on it. 
However there are many properties of the polar moment of inertia that enter the proof.
First it is an additive domain functional,
though only superadditivity is used. The first properties are
$$D(\Omega)+D(\Omega^*)\ge D(\Omega\cup\Omega^*)+ D(\Omega\cap\Omega^*)
\qquad{\rm and}\qquad
D(\Omega)\ge D(\Omega^*) \ \forall\Omega\ 
 ,
$$
and the final property is the specific form of $I(\Omega^*)$ in terms of area. 
(Some slight generalization might be possible, e.g. for higher-order moments
$$ I_p(\Omega,z)=\int_\Omega |z- z'|^p , \qquad  
I_p(\Omega^*)=\frac{2\pi}{p+2} \left(\frac{|\Omega|}{\pi}\right)^{(p+2)/2} ,
$$
satisfies the properties mentioned in general for domain functional $D$
but a more elaborate relation involving $|\Omega\setminus\Omega^*|$ arises at the
end of the calculation.)

\medskip

We record the items below as other inequalities can be found.
Cauchy-Schwarz gives the rightmost of
$$ \frac{(|\Omega|-|\Omega\setminus{\Omega^*_z}|)^2}{2\pi}
=\frac{|\Omega\cap{\Omega^*_z}|^2}{2\pi}
\le  \int_{\Omega\cap{\Omega^*_z}}  |z-z'|^2 
= \int  |z-z'| \chi_\Omega\ |z-z'| \chi_{\Omega^*_z}
\le \sqrt{ I_z(\Omega) \frac{|\Omega|^2}{2\pi} }, 
$$ 
where the characteristic functions $\chi$ have arguments $z'$.
Making use of the fact that the
square of the radius of ${\Omega^*_z}$ is $|\Omega|/\pi$, 
\begin{align*}
\alpha(\Omega)^2 \frac{|\Omega|^2}{2\pi}
\le \frac{|{\Omega^*_z\setminus{\Omega} }|^2}{2\pi}
&\le \int_{\Omega^*_z\setminus{\Omega} } |z-z'|^2 
\le \frac{|{\Omega^*_z\setminus{\Omega} }| |\Omega|}{\pi} ,
\\
 2 \alpha(\Omega) \frac{|\Omega|^2}{2\pi}
\le \frac{|\Omega\setminus{\Omega^*_z} | |\Omega|}{\pi} 
 &\le  \int_{\Omega\setminus{\Omega^*_z} } |z-z'|^2 
\end{align*}
 

\medskip\par\noindent{\bf An aside.}
Let $\overline{r}(\Omega)$ denote the transfinute diameter (also known as
logarithmic capacity or outer mapping radius).
It is possible to use this, in tandem with the inequality
$$ I_c\le \frac{\pi}{2} {\overline{r}}^4 \qquad{\rm equivalently}\ \ 
\left(\frac{ \overline{r(\Omega)} }
            { \overline{r(\Omega^*)}  }\right)^4 
\ge \frac{I_c(\Omega)}{ I_c(\Omega^*)}.
$$
See~\cite{PoS51} p126.  
(Note the different uses of the overline here to later subsections,
e.g. the table.)\\
Rewriting an earler inequality relating $I_c$ to the perimeter
$|\partial\Omega|$,
$$ I_c 
\le \frac{\pi}{2} \, 
\left(\frac{|\partial\Omega|}{2\pi} \right)^4 \qquad{\rm equivalently}\ \
\left(\frac{ |\partial\Omega|}{\sqrt{4\pi |\Omega|}}\right)^4
\ge \frac{I_c(\Omega)}{ I_c(\Omega^*)} .
$$
Hence the classical geometric inequality, concerning the isoperimetric
deficit, is improved to
\begin{equation}
\frac{ |\partial\Omega|^2}{4\pi |\Omega|} - 1
\ge \sqrt{1+2\alpha(\Omega)^2} -1 .
\label{eq:deficitMoA}
\end{equation}

\medskip\par\noindent{\bf Regular polygons.}
(Amongst the domains considered in detail in this paper are equilateral triangles and squares.
See also \S\ref{RegularngonsBetaGe0}.)
The polar moment of inertia about the centre of a regular $n$-gon is
$$ I_c(n{\rm -gon})
= \frac{|\Omega|^2}{6 n}\left(\tan(\frac{\pi}{n})+ 3\cot(\frac{\pi}{n})\right) . $$
In the literature, e.g.~\cite{PoS51} (where $n=3$, $4$ and $6$ is given on pages 252, 256 and 258),
the functionals are given in terms of side length, there denoted by $a$.
Then
$$ |\Omega| = \frac{n a^2}{4\tan(\frac{\pi}{n}) } . $$
The isoperimetric deficit satisfies
$${ |\partial\Omega|^2}-{4\pi |\Omega|}
= a^2 n^2\left( 1- \frac{\pi/n}{\tan(\pi/n)}\right) 
= 4|\Omega|\pi \frac{\tan(\pi/n)}{\pi/n} \left( 1- \frac{\pi/n}{\tan(\pi/n)}\right) 
\sim \frac{4|\Omega|\pi}{3} \left(\frac{\pi}{n}\right)^2 ,
 $$
 as $n$ tends to infinity.
 The circumradius $R$ is related to the area by
 $$|\Omega|=\frac{n R^2}{2}\sin(\frac{2\pi}{n}) . $$
 The inradius $\rho$ satisfies (by a simple application of Pythagoras'  Theorem)
 $$ R^2-\rho^2 = \frac{a^2}{4} , $$
 and hence is given by
 $$ \rho =\sqrt{R^2 -\frac{|\Omega|\tan(\pi/n)}{n}} .$$
 
 The modulus of asymmetry has been calculated for various polygons in~\cite{Fi14}.
 As an example
 $$2\alpha({\rm square})
 =\frac{4}{\pi} \left( 2\arccos(\frac{\sqrt{\pi}}{2}) -\frac{1}{2}\sqrt{\pi(4-\pi)}\right)
 \approx 0.1810919377 \ .
 $$
 Below we use the values of $\alpha(\Omega)$ from~\cite{Fi14} to compare the left-hand and right-hand sides
 of the inequality~(\ref{eq:deficitMoA}).
 \medskip
 \begin{center}
 \begin{tabular}{|| c | c | c | c | c ||}
\hline
$n$& $2\alpha(\Omega)$ & $\frac{ |\partial\Omega|^2}{4\pi |\Omega|} - 1$&
$\frac{\sqrt{2\pi I_c(\Omega)}}{|\Omega|} -1$& $\sqrt{1+2\alpha(\Omega)^2} -1 $\\
\hline
3& 0.3649426110&$\frac{3\sqrt{3}}{\pi}-1=0.653986686$& 0.099636111 &0.032759195\\
4& 0.1810919377&$\frac{4}{\pi}-1=0.273239544$&0.023326708  &0.008165237\\
6& 0.0744657545&$\frac{2\sqrt{3}}{\pi}-1= 0.102657791$&0.003825838 &0.001385327\\
 & & & &  \\
 $n$& & $\frac{\tan(\pi/n)}{\pi/n}-1$& $\sqrt{\frac{\pi}{3n}\left(\tan(\frac{\pi}{n})+ 3\cot(\frac{\pi}{n})\right)} -1$&  \\
$n>>1$&  &$\frac{\pi^2}{3 n^2}$ &$\frac{2\pi^4}{45 n^4}$ &   \\
$\infty$ &0 &0& 0 &0 \\
\hline
\end{tabular}
\end{center}

\subsubsection{Improving other isoperimetric inequalities in terms of $\alpha$}

When $\slipparameter=0$, the modulus of asymmetry 
$\alpha$ is used to bound the departure of $\lambda_1$ in~\cite{FMP09}.
In the inequality following their equation (2), for our situation of the ordinary Laplacian and the plane,
$$ \lambda_1(\Omega) \ge  \lambda_1(\Omega^*) \left( 1 +c\, \alpha(\Omega)^4 \right) $$
for some constant $c>0$.
(When estimating $c$, recall that our modulus of asymmetry differs by a factor of $2$ from that in~\cite{FMP09}.)

The isoperimetric deficit is used to
bound the departure of the eigenvalues in~\cite{Ni14}.
See also~\cite{MZ15}.
\medskip

\par\noindent{\bf Regular polygons, and $\beta=0$.} Data on these is available in~\cite{PoS51}, e.g. for $\lambda_1$:
\begin{itemize}
\item equilateral triangle $\lambda_1 |\Omega| = \frac{4\pi^2}{\sqrt{3}}=22.793$
\item square $\lambda_1 |\Omega| =2\pi^2=19.7392$
\item regular hexagon
\item
\item circle $\lambda_1 |\Omega| =\pi \, j^2 =  18.168$
\end{itemize}

\subsection{$\beta=0$ and bounds in terms of the classical isoperimetric deficit 
$|\partial\Omega|^2-4\pi |\Omega|$}
\label{subsec:deficit}

Defining $\Psi= 1 - 4\pi |\Omega|/|\partial\Omega|^2$, \cite{Pa67} inequality~(4.12)
gives 
$$ Q_{\rm steady}(\Omega)
\ge \frac{|\Omega|^2}{8 \pi} \left(
1- 2\Psi^2 (1-\Psi^2)^{-1} - 4\Psi^4 (1-\Psi^2)^{-2}\log(\Psi)
\right)\qquad
\mbox{\rm when\ } \beta=0 .
$$

There are also results for $\lambda_1$, and conjectures, e.g. the following.\\

{\par\noindent}{\bf Conjecture.} (Taken from~\cite{FK08}:see their equation~(12).)
{\it It is conjectured (in~\cite{AF06}) that for simply-connected plane domains
$$ \lambda_1(\Omega)
\le \frac{\pi j^2}{|\Omega|} +
\frac{\pi^2}{4} \, \frac{|\partial\Omega|^2 -4\pi|\Omega|}{|\Omega|^2} .
$$}\\
Equality occurs not only for the disc but also asymptotically on infinite rectangular strips.

\subsubsection{Relations between the classical isoperimetric deficit and $\alpha(\Omega)$}

See~\cite{FMP08}.

\medskip

\section{Nearly circular domains}\label{sec:nearCirc}

\subsection{Geometric preliminaries, nearly circular domains}
\label{subsec:nearCirc}

The domain functionals $Q_{\rm steady}$ and $\lambda_1$ involve pde problems,
and much effort has been expended seeking to bound these in terms of
simpler geometric functionals.
Thus Theorems \ref{thm:thm1} and \ref{thm:thm2} give inequalities
involving the area, $|\Omega|$.
Besides area, there are other geometric functionals, perimeter,
polar moment of inertia about the centroid, etc.

In much of this section we treat boundaries whose equation in polar form is
$$ r = a (1+ \rho(\theta) )\qquad{\rm with \ }a=1  .$$
For most of this we consider $\rho$ small.
In connection with small perturbations for a circle,
we follow the methods given in~\cite{PoS51}, p53.
There the boundary of the domain is given by
$$ r= 1+ \delta {\overline\rho}(\theta)
$$
with $\delta$ a first-order small quantity.
The boundary function is given as a Fourier series
\begin{equation}
\rho(\theta)
= a_0 + \sum_{n=1}^\infty (a_n\cos(n\theta) + b_n \sin(n\theta) ) .
\label{eq:rhoFourier}
\end{equation}
(Caution. The Fourier coefficients with $n\ge{1}$  in~\cite{PoS51} equation (1) are half ours as they have
a factor of 2 outside their sum.)

Various geometric functionals are considered in~\cite{PoS51}:
perimeter $|\partial\Omega|$, polar moment of inertia $I_c$, etc..
These can be approximated using $\rho$ small:
We begin with considerations of the area.

\subsubsection{Area $|\Omega|$}

With the boundary described by equation~(\ref{eq:rhoFourier}) the area of $\Omega$ is given by
\begin{equation}
\frac{|\Omega|}{\pi}
=  (1+a_0)^2 +\frac{1}{2}  \sum_{n=1}^\infty  (a_n^2+ b_n^2) .
\label{eq:areaFourier}
\end{equation}
\medskip


\medskip
\par\noindent{\bf Ellipse}
\smallskip

A test case for the methods is the ellipse $\frac{x^2}{a^2}+a^2 y^2 \le{1}$.
In polar coordinates relative to the centre the boundary curve is
\begin{align}
 r
&= \frac{1}{\sqrt{\frac{\cos(\theta)^2}{a^2} + a^2 \sin(\theta)^2}}
= \sqrt{\frac{2}{{a}^{2}+{a}^{-2}-\left( {a}^{2}-{a}^{-2} \right) {\cos(2\theta)}}} 
\nonumber \\
&= r\left(\frac{\pi}{4}\right) 
\left( 1 - 
\frac{a^2-a^{-2}}{a^2+a^{-2}} 
\cos(2\theta)
\right)^{-1/2}\ \ 
{\rm where}\ \  r\left(\frac{\pi}{4}\right) = \sqrt{\frac{2}{{a}^{2}+{a}^{-2}}}
\label{eq:rEllipsec2} \\
&=  \sqrt{\frac{1 +\tan^2(\theta)}{a^{-2} + a^2 \tan^2(\theta)} }
\label{eq:rEllipset} 
\end{align}
Without any assumptions on $a$, the coefficients in the Fourier series for $r(\theta)$ involve elliptic integrals. 
Our interest is in $a$ near 1.
However $(a-1)$ might not be the best perturbation parameter.
We have also used, in computations,
$$\epsilon
=\frac{a^2-a^{-2}}{a^2+a^{-2}} 
\sim 2(a-1) \qquad {\rm\ as\  }\ a\rightarrow{1} ,
$$
and others, e.g.~\cite{Day55}, use the eccentricity
$$ e=\sqrt{1-\frac{1}{a^4}} \sim 2(a-1)^{1/2}  \qquad {\rm\ as\  }\ a\rightarrow{1} .
$$
Including an $a$ dependence in $r(a,\theta)$ we remark that $r(1/a,\theta)=r(a,\theta+\pi/2)$ and
$\epsilon(1/a)=-\epsilon(a)$.
The binomial expansions
\begin{align}
\frac{r(0)}{r(\pi/4)}
&= (1-\epsilon)^{-1/2}
= 1+ \sum_{k=1}^\infty \frac{ (2k)! }{2^{2k} (k!)^2}\ \epsilon^k , 
\nonumber \\
\frac{r(\theta)}{r(\pi/4)}
&= (1-\epsilon\cos(2\theta))^{-1/2}
= 1+ \sum_{k=1}^\infty \frac{ (2k)! }{2^{2k} (k!)^2}\ \cos(2\theta))^k \epsilon^k , 
\label{eq:rCosk}
\end{align}
may be useful in finding higher terms in the perturbation expansions of some domain functionals.
The symmetries of the ellipse explain the form of the expansion:
\begin{itemize}
\item It is symmetric about $\theta=0$, hence only cosine terms.
\item It is symmetric about $\theta=\pi/2$ and hence only the even order cosine terms  $\cos(2 m\theta)$.
\item When $\epsilon$ is replaced by $-\epsilon$ and $\theta$ by $\theta+\pi/2$ the expression is unchanged and hence the form of the polynomial coefficients in $\epsilon$ forming the Fourier coefficients. 
For $m$ is odd, only odd powers of $\epsilon$ appear: for $m$ is even, only even powers of $\epsilon$ appear.
\end{itemize}
We will see these symmetries in connection with $\phi_1$ and $u_{\rm steady}$.
Returning to the study of the boundary curve,
we also need the expansion for $r(\pi/4)$:
$$ a^2=\sqrt{\frac{1+\epsilon}{1-\epsilon}} , \qquad  
\frac{2}{r(\pi/4)}= a^2 + a^{-2} =\frac{2}{\sqrt{1-\epsilon^2}},\qquad
r\left(\frac{\pi}{4}\right) = \left(1 -\epsilon^2 \right)^{1/4} .
$$
The first few terms give
$$ r(\theta)
= 1+\frac{1}{2}\,\cos \left( 2\,\theta \right) \epsilon+ \left( \frac{3}{16}\,\cos \left( 4\,\theta \right) -\frac{1}{16} \right) {\epsilon}^{2} +
 O(\epsilon^3) .
$$
Alternatively we can consider asymptotics as $a\rightarrow{1}$. 
Then $\rho=r_{\rm ellipse}-1$ satisfies 
{\small
\begin{align}
\rho(\theta) 
&\sim (a-1) \cos(2\theta) - (\frac{1}{4} +\frac{1}{2} \cos(2\theta)-\frac{3}{4} \cos(4\theta))  (a-1)^2 + O(  (a-1)^3) 
\nonumber \\
&= -\frac{1}{4} (a-1)^2 +  \left((a-1)-\frac{1}{2}  (a-1)^2\right) \cos(2\theta)+\frac{3}{4}  (a-1)^2\cos(4\theta)  + O(  (a-1)^3)   .
\label{eq:ellipseRho}
\end{align}
}
in the notation of equation~(\ref{eq:rhoFourier})
$$ a_0 =  -\frac{1}{4} (a-1)^2  +   o(  (a-1)^2) , \qquad
a_2 =  (a-1) + o( (a-1)) ,
$$
and all other Fourier coefficients are of sufficiently small order that they are not needed in subsequent calculations in this paper.
Of course we already know that the area of our ellipse is $\pi$.
The Fourier coefficients given in equation~(\ref{eq:ellipseRho}) for the approximation to the ellipse's boundary 
inserted  into equation~(\ref{eq:areaFourier}) is consistent with this:
\begin{align*}
\frac{\rm areaEllipse}{\pi}
&\sim \left(1- \frac{1}{4}(a-1)^2\right)^2 + \frac{1}{2}( a-1)^2 +O( (a-1)^3) , \\
&= 1  +O( (a-1)^3) .
\end{align*}

\medskip
For calculations extending the use of higher order terms in the Fourier series~(\ref{eq:ellipseRho}).
one might need formulae like
 $$   x^{2n}=2^{1-2n}\left(\frac{1}{2} \binom{2n}{n}+\sum_{j=1}^n \binom{2n}{n-j}T_{2j}(x) \right), 
 $$
and we plan, for elliptical $\Omega$ to implement this but have yet to do so.

\medskip
We now consider general curves close to circular with the $a_0$, $a_n$ and $b_n$ all small, $O(\delta)$ or smaller.
Equation~(\ref{eq:areaFourier}) is
$$
\frac{|\Omega|}{\pi}
=  (1+a_0)^2 \left( 1 +\frac{\sum_{n=1}^\infty  (a_n^2+ b_n^2) }{2 (1+a_0)^2}\right) ,
$$
and hence, with the $\dots$ indicating smaller order terms
\begin{align}
\sqrt{\frac{|\Omega|}{\pi}}
&\sim  (1+a_0) \left( 1 + \frac{\sum_{n=1}^\infty  (a_n^2+ b_n^2) }{4 (1+a_0)^2} + \dots \right) ,
\nonumber\\
&\sim (1+a_0)  +\frac{1}{4}  \sum_{n=1}^\infty   (a_n^2+ b_n^2)  + \dots \ .
\label{eq:areaGen}
\end{align}

\medskip
\subsubsection{Perimeter $|\partial\Omega|$}

The approximation for the perimeter when $\rho$ is small is as follows.
(Note the Caution given earlier about the difference in notation of our Fourier coefficients and those in~\cite{PoS51}.)
\begin{align}
|\partial\Omega|
&= \int_0^{2\pi} \sqrt{ r^2 + (\frac{d r}{d\theta})^2} \,  d\theta
= \int_0^{2\pi} \sqrt{ 1+2\rho +\rho^2 + (\frac{d\rho}{d\theta})^2} \  d\theta
\nonumber \\
&\sim 2\pi \int_0^{2\pi} \left( 1+\frac{1}{2} (2\rho +\rho^2 + (\frac{d\rho}{d\theta})^2) - \frac{1}{8} (2\rho)^2 \right) \ d\theta
\nonumber  \\
&\sim 2\pi \int_0^{2\pi} \left( 1+\rho  + \frac{1}{2}(\frac{d\rho}{d\theta})^2)  \right)\ d\theta
\nonumber  \\
\frac{|\partial\Omega|}{2\pi}
&\sim 1 + a_0 + \frac{1}{4}\sum_{n=1}^\infty  n^2\left( a_n^2 + b_n^2\right) .
\label{eq:perimGen}
\end{align}
\medskip

\par\noindent{\bf Ellipse}

The perimeter of our ellipse, area $\pi$, can be given without any approximation of small eccentricity $e$, by
$$|\partial\Omega|
= 4 a\,  {\rm EllipticE}(e)
\sim 2\pi a \left(1 -\frac{e^2}{4} - \frac{3 e^4}{64} \ \ldots\ \ \right)\ \qquad
{\rm for\ }\ \ a \rightarrow 1 . $$
Thus with $e^2\sim{2\epsilon}-2\epsilon^2$ and $a\sim{1+\epsilon/2+{\epsilon^2}/8}$,
\begin{equation}
|\partial\Omega|
\sim 2\pi(1+ \frac{3}{16}\epsilon^2 +\ \ \ldots ) \qquad
{\rm for\ }\ \epsilon \rightarrow 0 .
\label{eq:ellipsePerimEps}
\end{equation}

\medskip
\subsubsection{Polar moment of inertia about the centroid $I_c$}

The polar moment of inertia is, taking the origin at the centroid is
$$ I_c = \int_\Omega \left( (x-x_c)^2 + (y-y_c)^2 \right) ,
$$
where $z_c=(x_c,y_c)$ is the centroid of $\Omega$.
When the boundaries are given in polar coordinates, this is
$$I_c =\frac{1}{4}  \int_0^{2\pi} r(\theta)^4\ d\theta .
$$
\medskip\par\noindent{\bf Ellipse}

For our disk and ellipse these are
\begin{equation}
I_c({\rm disk})= \frac{\pi}{2} a^4, \qquad
I_c({\rm ellipse})= \frac{\pi}{4} (a^2 + a^{-2}) 
=\frac{\pi}{2\sqrt{1-\epsilon^2}} .
\label{eq:icEllipseEps}
\end{equation}
The asymptotics below check with the entries in the table of ~\cite{PoS51} reproduced in a later subsection:
\begin{align*}
\frac{2 I_c({\rm ellipse})}{\pi}
&\sim 1 + 2(a-1)^2 \qquad{\rm as\ } a\rightarrow 1 ,\\
\left(\frac{2 I_c({\rm ellipse})}{\pi}\right)^{1/4}
& \sim 1 + \frac{1}{2}(a-1)^2  + o(a-1)^2 ) \qquad{\rm as\ } a\rightarrow 1  ,\\
& \sim 1 + a_0 + \frac{3}{4}a_2^2    + o(a-1)^2 )\qquad{\rm as\ } a\rightarrow 1 .
\end{align*}

\subsection{Nearly circular domains with $\slipparameter=0$}
\label{subsec:nearCircBeta0}

First we treat the case  $\slipparameter=0$ and seek approximations,
notably for $Q_{\rm steady}$ and $\lambda_1$.
We follow~\cite{PoS51} in expressing the asymptotics (at small $\delta$) as
$${\mbox{\rm Quantity}} = 1+ a_0 + \frac{1}{4} \sum_{n=1}^\infty R(n) (a_n^2+ b_n^2) .$$
(The factor  $\frac{1}{4}$ is because of differing definitions of $a_n$ and $b_n$.)
Use an overline, as in~\cite{PoS51}, p33 to denote various radii. 

\subsubsection{The eigenvalue $\lambda_1$ at $\beta=0$}

\cite{PoS51} report that Rayleigh found
$$\frac{j}{\sqrt{\lambda_1}}
={\overline{\lambda_1}}
= 1 + a_0 -  \frac{1}{4} \sum_{n=1}^\infty \left( 1+ \frac{2 j {J_n}'(j)}{J_n(j)}\right) (a_n^2+ b_n^2) .
$$
Here the $\overline{\lambda_1}$ denotes the radius of the disk with principal eigenvalue $\lambda_1$.
Numeric values of the $n$-th term, involving the Bessel functions, are 
readily computed.
The $n=1$ term is $-1$. The rest  are positive. Furthermore as shown in~\cite{PoS51} p133 equation (2)
$$ 2 n-3   \le  \left( 1+ \frac{2 j {J_n}'(j)}{J_n(j)}\right) < 2n +1 , $$
with equality (in the left hand expression) only when $n=1$.\\
At large $n$,  we have $ \left( 1+ \frac{2 j {J_n}'(j)}{J_n(j)}\right) \sim 2n +1$.
The expressions $ \left( 1+ \frac{2 j {J_n}'(j)}{J_n(j)}\right)$ can, at each $n$,  be written as a rational function of $j$.

For an ellipse one finds
\begin{align}
\frac{j}{\sqrt{\lambda_1({\rm ellipse})}}
&\sim 1 -\frac{1}{4} (a-1)^2 -\frac{1}{4}  \left( 1+ \frac{2 j {J_2}'(j)}{J_2(j)}\right) (a-1)^2 + o( (a-1)^2) 
\nonumber \\
\frac{\lambda_1({\rm ellipse})}{j^2}
&\sim 1 +  \frac{1}{2} (a-1)^2 +\frac{1}{2}  \left( 1+ \frac{2 j {J_2}'(j)}{J_2(j)}\right) (a-1)^2 + o( (a-1)^2)
\nonumber  \\
&\sim 1 +    \left( 1+ \frac{ j {J_2}'(j)}{J_2(j)}\right) (a-1)^2 + o( (a-1)^2)
\nonumber  \\
&\sim 1 + \left(\frac{j^2}{2}-1\right) (a-1)^2 + o( (a-1)^2) ,
\label{eq:rayleigh}\\
&\sim 1 + 1.891592982\  (a-1)^2 + o( (a-1)^2) .
\nonumber 
\end{align}
\medskip

Without the special interest in nearly circular, eigenvalues for elliptical domains have been 
investigated before, with Mathieu functions in the eigenfunctions:
see ~\cite{Tr73,Day55,He49}.
In particular~\cite{Day55} equation (3.5)  reports that~\cite{He49} gave the approximation:
$$ 
\frac{\sqrt{\lambda}}{a} =	j 
\left(1- \frac{e}{4} - c_4  e^4 - c_6 e^6 - c_8  e^8 \right), 
$$
where $e$ is the eccentricity, defined above, and the numerical coefficients are
$$ c_4 = 0.034640 , \qquad c_6 =  0-010355 , \qquad c_8 = 0.004650 \ .
$$
This can be rearranged so that it determines $\lambda$
\begin{align*}
\lambda
&= a^2 j^2 \left(1- \frac{e}{4} - c_4  e^4 - c_6 e^6 - c_8  e^8 \right)^2 , \\
&\sim j^2 \left( 1 + (3-32 c_4)\  (a-1)^2 + \ldots \qquad\right)\qquad
{\rm\ as\ } a\rightarrow{1} ,\\
&\sim j^2 \left( 1 +  1.89152\  (a-1)^2 + \ldots \qquad\right)\qquad
{\rm\ as\ } a\rightarrow{1}
\end{align*}
where we have just indicated the first nontrivial term in the approximation for $a\rightarrow{1}$
and note that, except that the final decimal place is wrong, this agrees with Rayleigh's
result, equation~(\ref{eq:rayleigh}).


\bigskip
An inequality given in~\cite{PoS51} p99 is
$$ j^2\le\lambda_1 \le \frac{j^2}{2} (a^2 +a^{-2})
\ \left( \sim j^2 (1 + 2(a-1)^2 +\ldots ) \right)  , $$
the asymptotics in the parentheses being for $a$ tending to 1.
The right hand side isn't asymptotically tight for $a\rightarrow{1}$ as the $2$ there
differs from the $1.89\ldots$ of Rayleigh's result~equation(\ref{eq:rayleigh}).

\subsubsection{$Q_{\rm steady}$ at $\beta=0$}

Now, at $\beta=0$,
\begin{align}
{\frac{8 Q_{\rm steady}(\rm ellipse)}{\pi}} 
&= \frac{2}{(a^2+a^{-2})}  = \sqrt{1-\epsilon^2} ,
\nonumber \\
&\sim 1 - 2(a-1)^2+ O( (a-1)^3) 
\quad
{\rm as\ }\ a\rightarrow 1  ,
\nonumber \\
&\sim 1 - \frac{1}{2}\epsilon^2+ O( \epsilon^3) 
\quad
{\rm as\ }\ \epsilon\rightarrow 0,
\label{eq:PexactEps} \\
\left({\frac{8 Q_{\rm steady}(\rm ellipse)}{\pi}} \right)^{1/4}
&\sim 1 - \frac{1}{2} (a-1)^2+ O( (a-1)^3) 
\quad
{\rm as\ }\ a\rightarrow 1 .
\label{eq:PexactAsy}
\end{align}
Applying the formula from~\cite{PoS51} with the ellipse's values, just $a_0$ and $a_2$ to the order treated,
\begin{align*}
\left({\frac{8 Q_{\rm steady}(\rm ellipse)}{\pi}} \right)^{1/4}
&\sim 1+a_0 - \frac{(4-3)}{4} a_2^2 \\
&\sim 1 - \frac{1}{4} (a-1)^2 -  \frac{1}{4}  (a-1)^2  
\quad
{\rm as\ }\ a\rightarrow 1 \\
&\sim 1 - \frac{1}{2} (a-1)^2   
\quad
{\rm as\ }\ a\rightarrow 1 . .
\end{align*}

\subsubsection{Other domain functionals, mostly geometric}

\bigskip

\begin{tabular}{|| c | c | c | l ||}
\hline
Quantity& $R(n)$ & 
Comments \\
\hline
& & \\
$\frac{|\partial\Omega|}{2\pi}$& $n^2$&see equation~(\ref{eq:perimGen})\\
& & \\
$\overline{r}$& $2n-1$&outer mapping radius\\
& & transfinite diameter\\
$\overline{I_c}=\left(\frac{2I_c}{\pi}\right)^{1/4}$& 
$3$& \\
& & \\
$\overline{|\Omega|}=\sqrt{\frac{|\Omega|}{\pi}}$& $1$& see equation~(\ref{eq:areaGen})\\
& & \\
$\overline{Q_{\rm steady}}= \left({\frac{8 Q_{\rm steady}}{\pi}}\right)^{1/4}$& $-2n+3$& $\Delta{p}=1$, $\beta=0$\\
 $=\overline{P}$ & &the $P$ of~\cite{PoS51} is $4 Q_{\rm steady}$\\
& &See~\cite{PoS51} p131\\
& & \\
$\overline{\lambda_1}=\frac{j}{\sqrt{\lambda_1}}$& 
$\left( -1-\frac{2 j{J'}_n(j)}{ J_n(j)} \right)$& due to Rayleigh, $\beta=0$\\
& &See~\cite{PoS51} p132\\
& & \\
$(\dot{r}\le)\  r_c$& $-2n-1$& $r_c$ is inner mapping radius\\
& &with respect to the centroid $z_c$ \\
\hline
\end{tabular}
\medskip
\par\noindent
The ordering in the table is that of the isoperimetric inequalities
$$ {\dot{r}}
\le \overline{\lambda_1}
\le 
\overline{Q_{\rm steady}}
\le \overline{|\Omega|}=\sqrt{\frac{|\Omega|}{\pi}}
\le \overline{I_c}
\le \overline{r}
\le \frac{|\partial\Omega|}{2\pi} .
$$
See~\cite{PoS51} supplemented by equation~(\ref{eq:KJ}).

\medskip

\subsection{Nearly circular ellipses with $\slipparameter>0$}
\label{subsec:nearCircBetaNe0}

The symmetries of the ellipse noted in \S\ref{subsec:nearCirc} are useful to reduce the amount of computation.
(Of course assuming more parameters does not stop the calculations proceeding:
it merely slows them as one discovers that the unnecessary parameter is zero.)
The solution $u_{\rm steady}$ is of the form:
$$ \frac{u_{\rm steady}}{\Delta{p}}
=  \frac{1}{4}\left(1-r^2\right) +\frac{\slipparameter}{2} 
+\sum_{k=1}^\infty s_k(\epsilon) r^{2k} \cos(2k\theta) .
$$
Furthermore on integrating this over the ellipse
\begin{equation}
\frac{Q_{\rm steady}(\slipparameter)}{\Delta{p}}
=  \frac{\pi}{8} (1+4\slipparameter) 
+  \sum_{k=1}^\infty q_k\epsilon^{2k}. 
\label{eq:ellipseQinf}
\end{equation}
The solution $(\lambda_1,\phi_1)$ is, with $\lambda_1=\gamma^2$,  of the form
\begin{align*}
\gamma
&= \gamma_0 +  \sum_{k=1}^\infty \gamma_k\epsilon^{2k} ,\\
\phi_1
&= J_0(r\gamma)+  \sum_{k=1}^\infty a_k(\epsilon) J_{2k}(r\gamma) \cos(2k\theta) ,
\end{align*}
where, of course,  $\gamma_0$ is, as given in~\S\ref{sec:circular},  the least positive root of
$$ J_0(\gamma_0)= \beta \gamma_0 J_1(\gamma_0) . $$
That $Q_{\rm steady}$ and $\gamma$ are even in $\epsilon$ is evident from the
symmetry that changing the sign of $\epsilon$ is the same as interchanging $a$ and $/1a$:
geometrically the ellipse is the same.
The cosine series in $u_{\rm steady}$ and $\phi_1$ is because of the symmetry across the $x-axis$,
i.e. about $\theta-0$.
Also the presence of just the even order Fourier cosine coefficients is because of the symmetry across
the $y-axis$, i.e. about $\theta=\pi/2$.
There is also information on the structure of the $a_k(\epsilon)$ and $s_k(\epsilon)$ 
associated with the invariance under the simultaneous transformations,
$\epsilon\rightarrow{-\epsilon}$, $\theta\rightarrow{\theta+ \pi/2}$.
This requires that, 
when $k$ is odd, the functions $a_k(\epsilon)$ and $s_k(\epsilon)$ are odd in $\epsilon$;
when $k$ is even, the functions $a_k(\epsilon)$ and $s_k(\epsilon)$ are even in $\epsilon$.

\medskip

\goodbreak

\subsubsection{Details for the steady solution, $\beta\ge{0}$}

We present the $\beta=0$ solution here in order to check against our asymptotics for
the solution when $\beta\ge{0}$.
The $\beta=0$ solution dates back at least to St Venant.
In polar coordinates with the origin at the centroid of the ellipse
$$ u 
=\frac{1}{4}\left(   \sqrt {1-\epsilon^{2}} -r^2 +\epsilon \, r^2 \, \cos(2\theta)\, \right) ,
$$
and
$$Q_{\rm steady}({\rm ellipse}, \beta=0)
=\frac{\pi}{4 (a^2 + a^{-2})}
= \frac{\pi}{8} \, \sqrt {1-\epsilon^{2}} .
$$
When $\epsilon\rightarrow{0}$, this agrees with the asymptotics given earlier in
equation~(\ref{eq:PexactEps}).

Concerning the ellipse when $\beta>0$, the steady flow is not yet available, exactly, in general. 
Approximations -- some involving series and Mathieu functions -- are given in~\cite{DM07,DS16}.
When $\epsilon$ is small,  $u$ can be approximated in the form
$$ u =  \frac{1}{4}\left(  1-r^2 \right) + \frac{\beta}{2} +\epsilon^2\, t_{02} +
\epsilon \, t_{11} r^2 \cos(2\theta) +
\epsilon^2  \, t_{22} \, r^4\, \cos(4\theta) .
$$
The result of the perturbation analysis is
\begin{align*}
t_{11}
&= \frac{1}{4}\, \frac{1+\beta}{1+2\beta}
 , \\
t_{02}
&=  -\frac{1}{32}\, \frac{4 + 5\beta +6\beta^2}{1+2\beta}
=  -\frac{1}{32}\,\left( 1+3\beta + \frac{3}{1+2\beta}\right)
 , \\
 t_{22}
&=   -\frac{1}{32}\, \frac{\beta (1-2\beta)}{(1+4\beta)(1+2\beta)}
 .
\end{align*}
Integrating this over the ellipse gives the coefficients of the
expansion in~(\ref{eq:ellipseQinf}) for $Q_{\rm steady}({\rm ellipse})$:
$$Q_{\rm steady}({\rm ellipse})
\sim \frac{\pi}{8} (1+4\beta) + q_1\epsilon^2\ \ {\rm where} \ \ 
 q_1 = - \frac{\pi}{16} \left(1 - 8 t_{11} - 16 t_{02}\right)
 =  - \frac{\pi}{16} \left(1 +{\frac {\beta\, \left( 1+6\,\beta \right) }{2(2\,\beta+1)}}\right)  .
 $$
 When $\beta={0}$, this agrees with the asymptotics given earlier in
equation~(\ref{eq:PexactEps}).
When $\beta$ is large, $Q_\steady$ is large, as expected.

A task which remains is to check our results against those in~\cite{DM07}.

\medskip

\subsubsection{Details for the principal eigenvalue $\lambda_1$ and eigenfunction $\phi_1$ when $\beta\ge{0}$}

We seek the asymptotic approximation of the form
\begin{align*}
\gamma
&= \gamma_0 +  \gamma_2\, \epsilon^{2} ,\\
\phi_1
&= J_0(r\gamma)+   a_{11}\, \epsilon J_{2}(r\gamma) \cos(2\theta) 
 +   a_{22}    \, \epsilon^2 J_{4}(r\gamma) \cos(4\theta) ,
\end{align*}
where $\gamma_0$ is known from the solution for a circular disk.
Substituting this into the boundary condition and considering the successive  terms of the Taylor series in $\epsilon$
first yields
$$ a_{11} 
=\frac{{{\gamma_0}}^{2}}{4} \,{\frac { \left( {\beta}^{2}{{\gamma_0}}^{2}-\beta+1\right) }
 {{\beta}^{2}{{\gamma_0}}^{2}-2\,\beta+1}} ,$$
then (from simultaneous linear equations) $[a_{22},\gamma_2]$:
\begin{align}
 a_{22}
 &= {\frac {{{\gamma_0}}^{4}}{128}}\,{\frac{n_{22}}{d_{22}}} \ ,
\nonumber\\
 n_{22}
 &=  {{\gamma_0}}^{4} \left( {{\gamma_0}}^{2}-12 \right) {\beta}^{4}-{{\gamma_0
}}^{2} \left( 2\,{{\gamma_0}}^{2}-31 \right) {\beta}^{3}+ \left( -14-17\,{{
\gamma_0}}^{2}+2\,{{\gamma_0}}^{4} \right) {\beta}^{2}+ 
\nonumber\\
&\ \qquad\qquad\left( 25-2\,{{ \gamma_0}}^{2} \right) \beta+
{{\gamma_0}}^{2}-6 \ ,
\nonumber
\\
d_{22}
&=
\left( {\beta}^{2}{{\gamma_0}}^{2}-2\,\beta+1 \right) 
\left( -12\,{\beta}^{2}{{\gamma_0}}^{2}+{{\gamma_0}}^{4}{\beta}^{2}-2\,
\beta\,{{\gamma_0}}^{2}+24\,\beta+{{\gamma_0}}^{2}-6  \right) \ ,
\nonumber
\\
\gamma_2
&= \frac{ {\gamma_0}}{16}\,{\frac { \left(
{\beta}^{4}{{\gamma_0}}^{6}-{{\gamma_0}}^{2} \left( 2\,{{\gamma_0}}^{2}-3
 \right) {\beta}^{3}+ \left( 2\,{{\gamma_0}}^{2}+3 \right)  \left( -2+{{\gamma_0}}^{2} \right)
 {\beta}^{2}+ \left( 5-2\,{{\gamma_0}}^{2} \right) \beta-2+{{\gamma_0}}^{2}
\right)}
{ \left( {\beta}^{2}{{\gamma_0}}^
{2}-2\,\beta+1 \right)  \left( {\beta}^{2}{{\gamma_0}}^{2}+1 \right) }}\ .
\label{eq:gamma2beta}
\end{align}
 
A check on the solution is as follows.
Equation~(\ref{eq:rayleigh}) rewritten in terms of $\epsilon$ is
$$ \frac{\lambda_1({\rm ellipse})}{j^2}
= 1 + \frac{1}{8} (j^2-2) \epsilon^2 ,
$$
so that
$$ \frac{\sqrt{\lambda_1({\rm ellipse})}}{j}
= 1 + \frac{1}{16} (j^2-2) \epsilon^2 ,
$$
Setting $\beta=0$ in equation~(\ref{eq:gamma2beta}) gives
$$\frac{\gamma}{\gamma_0}
= 1 + \frac{\gamma_2}{\gamma_0}\epsilon^2
=1+  \frac{1}{16} (\gamma_0^2-2)\epsilon^2 ,
$$
which agrees with equation~(\ref{eq:rayleigh}) as when $\beta=0$, $\gamma_0=j$.

A further check on the solution is as follows.
When $\beta$ is large, an alternative approach begins with
$$\gamma^2=\lambda_1
\sim\frac{|\partial\Omega|}{\beta|\Omega|} \qquad
{\rm for\ }\ \beta\rightarrow\infty .
$$
Applying this formula to the nearly circular ellipses, with area $\pi$, on using
equation~(\ref{eq:ellipsePerimEps}), we have, in the double limit
$$\gamma^2
\sim\frac{2}{\beta}\, (1 +\frac{3}{16}\epsilon^2 ) ,\qquad
\frac{\gamma}{\gamma_0}
\sim (1 +\frac{3}{32}\epsilon^2 ) .
$$
On using $\gamma_0\sim\sqrt{2/\beta}$ in equation~(\ref{eq:gamma2beta}) we find
$$\frac{\gamma_2}{\gamma_0}
\sim
\frac{3}{32} ,
$$
which checks.

\section{Discussion and open questions}\label{sec:Discussion}

This paper reports  isoperimetric results found 
outside the context of microchannel flows.
We hope that
 the open questions of \S\ref{subsubsec:PolygonalbGt0} and more stated below may be of interest to --
 and resolved by readers of our paper.

The mathematical literature on this sort of pde problem has value for  engineers. 
Bounds for pde problems allow for checks on numeric computations.

For applied mathematicians there are many challenges.
It was noted in \S\ref{sec:Intro} that there is considerable uncertainty
about modelling of slip flows, and it may be that, under appropriate restrictions
on $\beta(u)$, it may be possible to prove Theorems 1 and 2 when the
boundary condition is
\begin{equation}
 u +\beta(u) \frac{\partial u}{\partial n} =0 .
\label{eq:betau}
\end{equation}
Another direction for generalization is replacing the Laplacian in the
pde with some other second-order elliptic operator.
Indeed some of the isoperimetric results have been established for the
$p$-Laplacian.
There may be non-Newtonian fluids (studied with different goals in~\cite{ME08}
for example) for which the results can be established.

There are various related issues.
One of the issues is to what extent the isoperimetric inequalities can be improved. 
For example if the cross-section is in some sense close to circular
by how much do the functionals $Q_{\steady}$ and $\lambda_1$ depart from
those of a circular cross-section of the same area.
There is also a tradition of perturbation expansions for nearly circular domains.
The Faber-Krahn inequality (the $\beta=0$ case of Theorem 2) was
conjectured by Rayleigh partly on the basis of information
concerning nearly circular domains.
This is reported in~\cite{PoS51}.

\newpage

\appendix
\section{Polygonal cross-sections}\label{sec:PolygonalIso}

\subsection{$n$-gons: fixing $n\ge{3}$, area and $\beta\ge{0}$}

\subsubsection{$\beta=0$}
\label{subsubsec:PolygonalbEq0}

The following results when $\beta=0$ are consequences of results in~\cite{PoS51}.
\begin{itemize}
\item
At a given pressure gradient, for all triangular channels with given area,
that which maximises the steady flow is equilateral.
Consider flows starting from rest developing from a constant imposed pressure gradient:
for all triangular channels with given cross-sectional area, that which
has the slowest approach to the steady flow is equilateral.
\item
At a given pressure gradient, for all quadrilateral channels with given area,
that which maximises the steady flow is square.
Consider flows starting from rest developing from a constant imposed pressure gradient:
for all quadrilaterial channels with given cross-sectional area, that which
has the slowest approach to the steady flow is square.
\end{itemize}
Steiner symmetrisation arguments can be used to prove these.
The first sentence is akin to the more familiar St Venant Inequality.
The second sentence is akin to the more familiar Faber-Krahn inequality.

See~\cite{PoS51} 
page vii and page 158 where this is set in a larger context.\\
Define {\it ${\cal{S}}(n,{\rm functional},\langle{\rm smallest}|{\rm largest}\rangle)$ as the statement:\\
Of all polygons with $n$ sides in some class $C$ and with a given area, the regular polygon (in $C$) has
the $\langle{\rm smallest}|{\rm largest}\rangle$ functional.}\\
A famous conjecture of Polya and Szego is that for any $n$\\
$\bullet$ ${\cal{S}}(n,{\rm functional},{\rm smallest})$ is true for the functional being any of
perimeter, moment of inertia, ... and fundamental Dirichlet eigenvalue $\lambda_1$
over the class $C$ of all $n$-gons, and\\
$\bullet$  ${\cal{S}}(n,{\rm functional},{\rm largest})$ is true for the functional being 
(any of inner mapping radius and) 
torsional rigidity ($Q_\steady$)
over the class $C$ of all $n$-gons.\\
The underlined result on~\cite{PoS51} page 158 states that both the above are true when $n=3$ (triangles)
and when $n=4$ (quadrilaterals), in each case with the pde functionals using $\beta=0$.

Further results concerning triangles are given in~\cite{Si08}.

Motivated by the unsolved problems when $\beta>0$, we also consider special cases of what has already been proved for general quadrilaterals.
Consider next  trapeziums (including parallelograms) whose parallel sides are distance $b$ apart and symmetrising about a perpendicular to these parallel sides. One has the following.\\
{\it Amongst all trapeziums whose parallel sides are distance $b$ apart and whose area is fixed, that which has the least fundamental Dirichlet eigenvalue $\lambda_1$ is the symmetric trapezium.}\\
It would be possibe to explore numerically whether there might be a similar result with Robin boundary conditions.
There are many other upper and lower bounds on the Dirichlet $\lambda_1$:
see for example~\cite{He66,He03,FS10}.

 \subsubsection{$\beta>0$}
 \label{subsubsec:PolygonalbGt0}

As noted before,
the exact solutions for equilateral triangles (see \S\ref{subsec:equilatQs}, \S\ref{subsec:equilatLambda1}) 
and rectangles (see \S\ref{subsec:Rectangular}) are available
for the Robin boundary condition case (Navier slip condition).
However symmetrisation techniques are inappropriate when $\beta>0$,
and we do not know if the results for triangles and quadrilaterals stated for $\beta=0$ in the preceding subsection
are also true when $\beta>0$.
\medskip

\par\noindent{\bf Triangles}\\
Indeed the question for triangles is noted as {\it Open Problem 1} in~\cite{LS15}.
Restated in our notation, the problem is:.\\
{\it Is ${\cal{S}}(3,{\lambda_1},{\rm smallest})$ true with $\lambda_1$ the fundamental Robin eigenvalue 
over the class $C$ of all triangles?}

\smallskip
A corresponding problem for $Q_\steady$ is also open.

\medskip

\par
\begin{figure}[hb]
\centerline{
$$
\includegraphics[height=5cm,width=6.5cm]{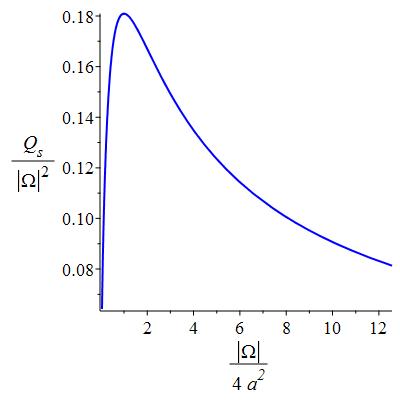}
\includegraphics[height=5cm,width=6.5cm]{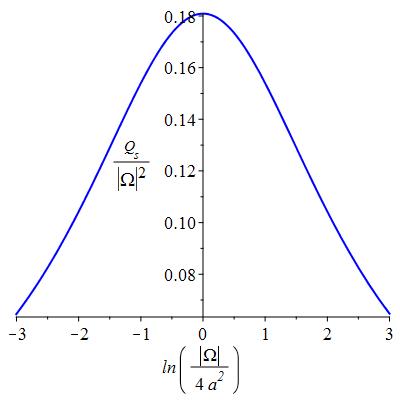}
$$
}
\caption{$\beta$  and $|\Omega|$ fixed. For rectangles with fixed area, the numerical evidence
(also noted in ~\cite{WWS14}) is that $Q_\steady$ is maximized by the square.
Shown are plots of $Q_\steady/|\Omega^2|$ (a) at left against $|\Omega|/(4 a^2)$
and (b) at right against $\log(|\Omega|/(4 a^2))$.}
\label{fig:QsRectb1Square}
\end{figure}

\begin{figure}[ht]
\centerline{
$$
\includegraphics[height=5cm,width=6.5cm]{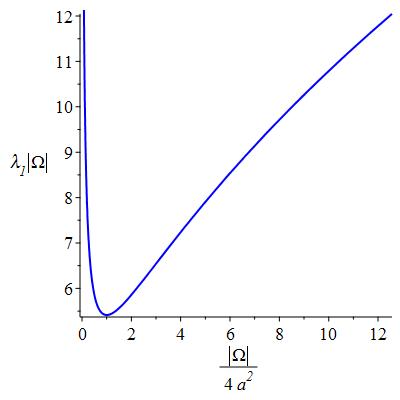}
\includegraphics[height=5cm,width=6.5cm]{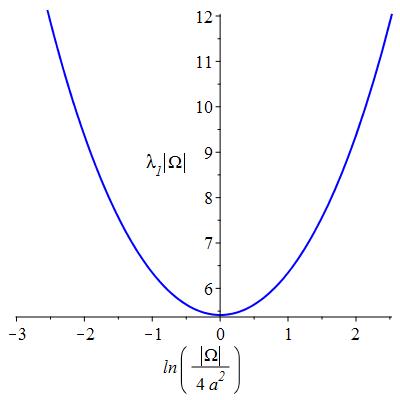}
$$
}
\caption{$\beta$. and $|\Omega|$ fixed.
 For rectangles with fixed area,  $\lambda_1$ is minimized by the square.
Shown are plots of $|\Omega|\lambda_1$ (a) at left,  against $|\Omega|/(4 a^2)$,
and (b) at right, against $\log(|\Omega|/(4 a^2))$}
\label{fig:lambda1Rectb1SquareMin}
\end{figure}

\noindent{\bf Rectangles}\\
The plots and result here are obtained from the formulae from \S\ref{subsec:Rectangular}.
In Figures~\ref{fig:QsRectb1Square} and \ref{fig:lambda1Rectb1SquareMin} we consider rectangles
of area $\pi$ as $a$ varies.
The numerics which led to these figures was independent of that in~\cite{WWS14}.
The Conclusion of~\cite{WWS14} presents the observation that there was numerical evidence 
that, over rectangles of a given area, $Q_\steady$ is maximized by the square.

\clearpage

The expressions for $Q_\steady$ are rather lengthy: the formulae for $\lambda_1$ are more tractable,
and the following small result, Theorem~\ref{thm:thm3},  consistent with the numerical evidence in 
the plot shown in Figure~\ref{fig:lambda1Rectb1SquareMin} is 
found by calculation.
 
\medskip

The main item for remainder of this paper is the proof of Theorem~\ref{thm:thm3}, which can be restated:\\
{\it
 ${\cal{S}}(4,{\lambda_1},{\rm smallest})$ is true with $\lambda_1$ the fundamental Robin eigenvalue 
over the class $C$ of all rectangles.}
Its proof is in \S\ref{sec:rectlambda1}.

\subsection{Regular $n$-gons: fixing  area and $\beta\ge{0}$, but varying $n$}
\label{RegularngonsBetaGe0}

In connection with $\beta=0$ and $\lambda_1$, see~\cite{Ni14a}.
The conjecture below is stated in~\cite{AF06}:\\
{\bf Conjecture.}  {\it
For all $n\ge{3}$ the first Dirichlet Laplacian eigenvalue of the regular $n$-gon is greater than the one of the regular $(n+1)$-gon of same area.}

\begin{figure}[ht]
\centerline{
$$
\includegraphics[height=5cm,width=6.5cm]{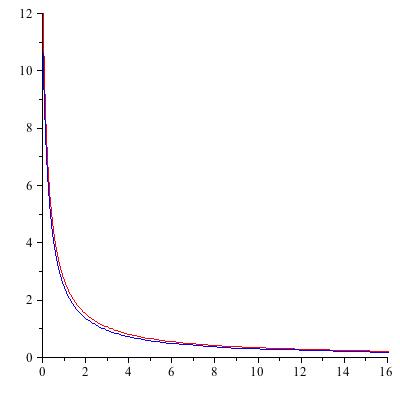}
$$
}
\caption{$\lambda_{1\Delta}$ (red) and $\lambda_{1\square}$ (blue) plotted against $\beta$.
Both areas are $\sqrt{3}$.}
\label{fig:lambda1TriSq}
\end{figure}

It is easy to verify this when $n=3$, as the formulae for both the equilateral triangle, and the square, are available.
See equations~(\ref{eq:betaLamTri}) and~(\ref{eq:betaLamSq}).
This suggests the question:\\
{\bf Question.}  {\it
Is it true that, for $\beta\ge{0}$, for all $n\ge{3}$ the first Robin Laplacian eigenvalue of the regular $n$-gon is greater than the one of the regular $(n+1)$-gon of same area.}\\
Again this is the case for $n=3$: see Figure~\ref{fig:lambda1TriSq}.

\smallskip

There is a corresponding question for $Q_\steady$. 
From the formulae we have for the equilateral triangle and the square:
$$ Q_{\steady\square}> Q_{\steady\Delta}, \qquad
\lambda_{1\square}<\lambda_{1\Delta} . $$



\clearpage

\section{Isoperimetric and other bounds on $\lambda_1$ for a rectangle}\label{sec:rectlambda1}

The plots and result here are obtained from the formulae from \S\ref{subsec:Rectangular}.
In 
Figure~\ref{fig:lambda1Rectb1SquareMin} we consider rectangles
of area $\pi$ as $a$ varies.

Our result, Theorem~\ref{thm:thm3},  is consistent with the numerical evidence in 
the plot shown in Figure~\ref{fig:lambda1Rectb1SquareMin}. 
\medskip

Much of the remainder of this section concerns the proof of Theorem~\ref{thm:thm3} and related inequalities.
We have several approaches to the proof. In all of them $\beta$ is given.\\
$\bullet$ The first of these in \S\ref{subsec:Proofexplicit}, is reasonably direct but uses a parametrisation which seems somewhat artificial.
(The geometry of the rectangle isn't give directly, but, $\mu_X $ and $\mu_Y$ are the given parameters.)
It is, so far, the only approach which has been carried through to complete the proof of Theorem~\ref{thm:thm3}.\\
$\bullet$  In \S\ref{subsec:variational} we begin a variational approach, which, at present, requires some numerics to
indicate that it will give a proof of the theorem.\\
$\bullet$  In a further approach, in \S\ref{subsec:hrSpecified},  the area is given, $4h^2$, as is the aspect ratio the rectangles, or at least $r$ which is
closely related to the aspect ratio.\\

 \smallskip
 \begin{theorem}{\label{thm:thm3}}{
Amongst all rectangles with a given area, that which has the smallest fundamental Robin eigenvalue
is the square.}
\end{theorem}

 \par\noindent{\it Comment.}
 The result, of course, is well known when $\beta=0$. Then
 $$ \lambda_1(r) =\frac{\pi^2}{4 h^2} \left( r^2 +\frac{1}{r^2} \right) \qquad
 {\mbox{\rm when \ \ }}  \beta=0 ,
 $$
 so
  $$ \frac{\lambda_1(r)}{ \lambda_1(1)}
  =\frac{1}{2} \left( r^2 +\frac{1}{r^2} \right) \qquad
 {\mbox{\rm when \ \ }}  \beta=0 .
$$
The function at the right of the equation above is minimized at $r=1$,
and monotonic either side of this.

One point of view is that we are required to establish that
$${\hat\mu}(r h)^2 + {\hat\mu}(\frac{h}{r})^2 \ge {\hat\mu}(h)^2 , $$
where the ${\hat\mu}$ function is only defined implicitly.
The inequality above might be rewritten, with $f(r)={\hat\mu}(r h)^2$ as
\begin{equation}
f(r) + f(\frac{1}{r}) \ge f(1) .
\label{eq:FI}
\end{equation}
It would be a valid, and potentially useful, exercise to find what functions
$f$, positive and decreasing on $0<r<\infty$, with$f((0,\infty))=(0,\infty)$,  satisfy such a (linear homgeneous) functional inequality, and then
show that $f(r)={\hat\mu}(r h)^2$ has properties so that it  is in the solution set of
(\ref{eq:FI}).
However an obstacle to this is the implicit definition of $\hat\mu$ and
our first proof avoids this by noting that equation~(\ref{eq:rectmuchat}) can be solved
explicitly for $\hat c$.

Our goal in the approach in \S\ref{subsec:hrSpecified} -- when $h$ and $r$ are prescribed --
is to show that the same properties can be established for the more elaborate 
expression~(\ref{eq:lambda1Rect}) applying when $\beta>0$.
Now, with the notation omitting dependence on $h$ and $\beta$, define
$$ E(r)=\frac{\lambda_1(r)}{ \lambda_1(1)}
  =\frac{ \left(\mu_1( r)^2 +\mu_1(\frac{1}{r})^2 \right)}{2\mu_1(1)^2} .
$$
There is some obvious symmetry, for example
$\lambda_1(r)=\lambda_1(1/r)$,
and this means it will suffice to establish that $E(r)\ge{1}$ for $0<r<1$.

Before treating the situation with $\beta$ general, we look at the asymptotics 
for $\beta$ small~(\ref{eq:lambda1Rectb0}) and
for $\beta$ large~(\ref{eq:lambda1RectbInf}).
For $\beta$ small the asymptotics in equation~(\ref{eq:lambda1Rectb0})
are clearly unsatisfactory as $r\rightarrow{0}$ or as $r\rightarrow{\infty}$
as the approximation to $\lambda_1$ so given becomes negative.
Thus the approximations require $\beta$ to be smaller order than $r$ or $1/r$.
The situation is more satisfactory with equation~(\ref{eq:lambda1RectbInf}).
Then
$$ E_\infty(r,\beta)\sim \frac{1}{2}\left(r +\frac{1}{r}\right)
=\frac{|\partial\Omega(r)|}{|\partial\Omega(1)|}  \qquad {\rm as\ \ }
\beta\rightarrow\infty
$$
and $E_\infty$ is clearly minimized at $r=1$.

\medskip
\goodbreak

\subsection{Proof of Theorem~\ref{thm:thm3} from explicit formulae in  $\mu_X$, $\mu_Y$}
\label{subsec:Proofexplicit}

In the approach in this subsection $\mu_X$ and $\mu_Y$ are specified, and these determine $h$ and $r$, and $\mu=\mu_\square$.

{
It will be useful to collect properties of a function occuring in one of the proofs of Theorem~\ref{thm:thm3}.
The function involves the $\arctan$ function and we note the following bounds on $\arctan(1/z)$
obtained by bounding the integrand below:
\begin{equation}
\frac{z}{1+z^2}=\frac{1}{z}\, \left(\frac{1}{1+\frac{1}{z^2}}\right)
<\int_0^{1/z}\frac{1}{1+t^2}\, dt< \frac{1}{z} 
\qquad {\rm for\ } z>0.
\label{eq:arctanIneq}
\end{equation}
The left hand inequality above is used in one method of proving part of Lemma 1.
\medskip
}{}

\par\noindent{\bf Lemma 1.} {\it Let
\begin{equation}
\phi_1(z)=\frac{1}{z} \, \arctan\left( \frac{1}{z} \right) ,\qquad
\phi_2(z)=\frac{1}{\sqrt{z}} \, \arctan\left( \frac{1}{\sqrt{z}} \right)  .
\label{eq:phiarctan}
\end{equation}
Then both $\phi_(z)$ and $\phi_2(z)$ are monotonic decreasing on $0<z<\infty$,
so, in particular,
$$ \phi_1({\hat\mu}_-)\ge\phi_1({\hat\mu}_+) \Longrightarrow \ {\hat\mu}_-\le {\hat\mu}_+ \qquad
{\rm and\ }\ \
 \phi_2({\hat\mu}_-^2)\ge\phi_2({\hat\mu}_+^2) \Longrightarrow \ {\hat\mu}_-^2\le {\hat\mu}_+^2 .
$$
Both $\phi_1$ and $\phi_2$ are  log-convex, and, even stronger, both are completely monotone.
}

\par\noindent{\it Proof.} (This is much longer than needed. Going directly to the inverse Laplace transforms as given below gives the complete monotonicity and thence everything else. However we leave it as first written.)
We have
$$
\phi_1' 
= -\frac{\phi_1}{z} - \frac{1}{z(1+z^2)} ,
$$
establishing the $\phi_1$ is decreasing on $0<z<\infty$. 
Continuing, we have
\begin{eqnarray*}
\phi_1''
&=& \frac{2\phi_1}{z^2} + \frac{2+4 z^2}{z^2(1+z^2)^2} ,\\
\phi_1^2 (\log(\phi_1))''
&=&  \frac{\phi_1^2}{z^2} + \frac{2\phi_1}{(1+z^2)^2} -\frac{1}{z^2 (1+z^2)^2} ,\\
&=&   \frac{1}{z^2} \left(\phi_1^2 -\frac{1}{(1+z^2)^2}\right) +  \frac{2\phi_1}{(1+z^2)^2} ,
\end{eqnarray*}
and the expression within the large parentheses in the preceding equation is
positive by the left hand side of inequality~(\ref{eq:arctanIneq}).
Hence $\phi_1$ is log-convex on $0<z<\infty$.

We also have
$$ \phi_1(z) = \int_0^\infty \exp(-z t) {\rm Si}(t)\, dt , $$
where {\rm Si} is the sine integral
$$ {\rm Si}(t) = \int_0^t \frac{\sin(\tau)}{\tau} \, d\tau . $$
Since ${\rm Si}(t)>0$ for $t>{0}$, $\phi_1(z)$ is completely monotone.
(Definitions and properties of completely monotone functions  are given in the Appendix.
In particular, we remark that all completely monotone functions are log-convex.)

Since $\phi_2(z)=\phi_1(\sqrt{z})$, i.e. $\phi_2$ is the composition of a completely monotone
function with a Bernstein function (a positive function whose first derivative is completely monotone),
$\phi_2$ is completely monotone. (See Theorem 2 of~\cite{MS01}.)
More directly $\phi_2$ is completely  monotone as
$$ \phi_2(z)
= \int_0^\infty \exp(-z t)\, \frac{1}{2} \sqrt{\frac{\pi}{t}} {\rm erf}(\sqrt{t})\, dt ,
$$
and the integrand in the expression above is positive.\qed

\medskip

\par\noindent{\bf Lemma 2.} {\it Let $\zeta:(0,\infty)\rightarrow{\mathbf R}$ be twice continuously differentiable,
decreasing and convex.
Then the separable convex programming problem, given $Z>0$,  to find $X$ and $Y$ 
which minimize $X+Y$ and satisfy the constraint
$$ \zeta(X) +\zeta(Y)=  2\zeta(Z) , $$
is solved by $X=Z=Y$.}

\par\noindent{\it Proof.} The convexity of $\zeta$ gives
$$ \zeta\left(\frac{X+Y}{2}\right)
\le \frac{\zeta(X)+\zeta(Y)}{2}=\zeta(Z) .$$
But $\zeta$ is decreasing so
$$\frac{X+Y}{2} \ge Z , \qquad X+Y \ge 2Z .$$
The lower bound of $2Z$ on the sum $X+Y$ is attained when
$X=Z$ and $Y=Z$.\qed

\medskip

\par\noindent{\it Proof of Theorem~\ref{thm:thm3}.} 
\smallskip

We have already noted that at fixed $\beta$ and $h$, $\mu(r h)$ is monotonic in $r$, and in particular,
$$ \mu_Y = \mu(\frac{h}{r}) \le \mu(h)=\sqrt{\frac{\lambda_1({\rm square})}{2}} \le \mu(r h) = \mu_X
\qquad
{\rm for\ }\ 0<r\le{1}. $$
Let us now suppose we are given $\beta$, $\mu_X$ and $\mu_Y$ with $\mu_X>\mu_Y$.
Where there can be no confusion about the value of $h$ we sometimes abbreviate:
$$ \mu_Y= \mu(\frac{h}{r})  {\ \rm as\ } \mu(\frac{1}{r}),\qquad
\mu(h)  {\ \rm as\ } \mu(r=1) {\ \rm or\  as\ } \mu_\square, \qquad
\mu_X= \mu(r h) \  {\ \rm as\ }\mu(r) .
$$
As in (\ref{eq:rectmuchat}) also write
$$ \beta\mu_Y = {\hat\mu}_Y,\qquad \mu_\square = {\hat\mu}_\square, \qquad  \beta\mu_X = {\hat\mu}_X .$$
The goal is to show
$$ \mu_X^2 + \mu_Y^2 \ge  2\mu(r=1)^2 = 2\mu_\square^2 = \lambda_1({\rm square}) ,$$
or equivalently
\begin{equation}
{\hat\mu}_X^2 + {\hat\mu}_Y^2 \ge 2 {\hat\mu}_\square^2 . 
\label{eq:goalB}
\end{equation}
Now, in the notation of equation~(\ref{eq:phiarctan})
$$
\frac{hr}{\beta} = \phi_1({\hat\mu}_X)= \phi_2({\hat\mu}_X^2) ,\quad
\frac{h}{\beta} = \phi_1({\hat\mu}_\square) = \phi_2({\hat\mu}_\square^2), \quad
\frac{h}{r\beta} = \phi_1({\hat\mu}_Y)= \phi_2({\hat\mu}_Y^2)  . $$
Now, with $\beta$ given,  $h$ is determined by the values of $\mu_X$ and $\mu_Y$:
$$ \left(\frac{h}{\beta}\right)^2 
= \phi_1({\hat\mu}_X)\phi_1({\hat\mu}_Y)= \phi_2({\hat\mu}_X^2)\phi_2({\hat\mu}_Y^2). $$
Also
$$ \left(\frac{h}{\beta}\right)^2 
= \phi_1({\hat\mu}_\square)^2 = \phi_2({\hat\mu}_\square^2)^2. $$
Eliminating $h/\beta$ between the two preceding equations we have
$$ \phi_1({\hat\mu}_X)\phi_1({\hat\mu}_Y)= \phi_1({\hat\mu}_\square )^2 ,\qquad
 \phi_2({\hat\mu}_X^2)\phi_2({\hat\mu}_Y^2)= \phi_2({\hat\mu}_\square^2 )^2 .$$
 
 We now apply Lemma 2 with $\zeta(z)=\log(\phi_2(z))$ and note that Lemma 1
 ensures that the conditions on $\zeta$ needed for Lemma 2 are satisfied.
 This establishes the result.\qed

\bigskip
 
 A slightly stronger inequality follows from considering $\phi_1$.
Since $\phi_1$ is log-convex
$$\phi_1\left(\frac{{\hat\mu}_X + {\hat\mu}_Y}{2}\right)^2 
< \phi_1({\hat\mu}_X)\phi_1({\hat\mu}_Y)= \phi_1({\hat\mu}_\square )^2 . $$

\begin{figure}[ht]
\centerline{
$$
\includegraphics[height=5cm,width=6.5cm]{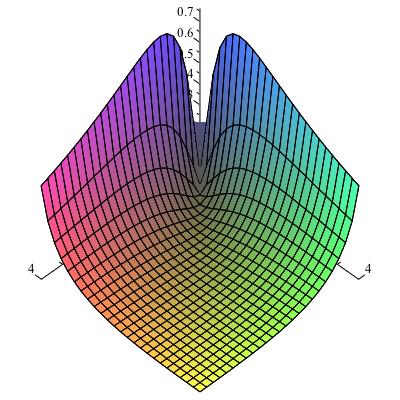}
$$
}
\caption{A plot of  $\sqrt{\phi_1({\hat\mu}_X)\phi_1({\hat\mu}_Y)}-\phi_1\left(\sqrt{\frac{{\hat\mu}_X^2 + {\hat\mu}_Y^2}{2}}\right)$ 
indicating that it is positive in the quadrant off the line ${\hat\mu}_X={\hat\mu}_Y$ where the expression is 0.}
\label{fig:logcvxB}
\end{figure}

A short, but reassuring,  digression at this point is to look at the plot given in
Figure~\ref{fig:logcvxB}.

Since $\phi_1$ is decreasing, the immediately preceding inequality shows
\begin{equation}
 \frac{{\hat\mu}_X + {\hat\mu}_Y}{2} >  {\hat\mu}_\square . 
\label{eq:logcvxB}
\end{equation} 
This is a stronger inequality than in our Theorem, and we deduce the Theorem (again) from
\begin{eqnarray*}
 \frac{{\hat\mu}_X^2 + {\hat\mu}_Y^2}{2} 
&=& \left(\frac{{\hat\mu}_X + {\hat\mu}_Y}{2} \right)^2 + \left(\frac{{\hat\mu}_X - {\hat\mu}_Y}{2} \right)^2 \\
&>&  \left(\frac{{\hat\mu}_X + {\hat\mu}_Y}{2} \right)^2 .
\end{eqnarray*}
Combining the last inequality with~(\ref{eq:logcvxB}) we have
$$\frac{{\hat\mu}_X^2 + {\hat\mu}_Y^2}{2} > \left(\frac{{\hat\mu}_X + {\hat\mu}_Y}{2} \right)^2> {\hat\mu}_\square^2 ,$$
which again reaches our goal~(\ref{eq:goalB}).\qquad$\square$

\subsection{A variational approach to the theorem and other inequalities}
\label{subsec:variational}

For domains $\Omega$ with sufficiently smooth boundaries (including our rectangles),
the fundamental eigenvalue is given by
$$\lambda_1
= {{\rm min}\atop{v\in{W^1_2}}}{\cal V}(v)\qquad{\rm where\ } \ \ 
{\cal V}(v)= \frac{G(v) + E(v)/\beta}{F(v)} ,
$$
and where the gradient term $G$ is
$$ G(v)  = \int_\Omega |\nabla v|^2,
$$
the face term $F$ is
$$ F(v)  = \int_\Omega v^2,
$$
and the edge term $E$ is
$$ E(v)  = \int_{\partial\Omega} v^2.
$$

In this subsection we will apply this to rectangles 
$\Omega= (-h r,h r)\times(-h/r,h/r)$ and with test functions
$v(x,y)=V_X(x)V_Y(y)$ and we establish inequalities of the form
\begin{alignat*}{3}
\lambda_\square
&= 2\mu_\square^2
&&\le f(\mu_X,r) + f(\mu_Y,\frac{1}{r}) ,\\
\lambda_1(r)
&= \mu_X^2 + \mu_Y^2
&&\le g(\mu_\square,r) + g(\mu_\square,\frac{1}{r}) .
\end{alignat*}
Here $f(\mu,1)=\mu^2$ and $g(\mu,1)=\mu^2$.
Often, as in the preceding sentence, $\mu$ is an arbitrary positive number.
However, where no misunderstanding is possible we will omit the subscript from
$\mu_\square$, writing it just as $\mu$.
The inequality of Theorem~\ref{thm:thm3} is an instance of the first type above, 
an upper bound on $\lambda_\square$, though we have yet to get, via the variational methods,
a neat proof of this. The more elaborate $f$ we find seem to be stronger inequalities.
As an example of the second type we mention the neat inequality~(\ref{eq:VxVymuSqSimple}).

For the rectangle, we write 
$$ G(v)= G_x(v) + G_y(v), \qquad
G_x(v) = \int_{-b}^b\int_{-a}^a v_y^2\, dx\, dy,\qquad
G_y(v) = \int_{-b}^b\int_{-a}^a  v_x^2,
$$
$$ F(v) = \int_{-b}^b\int_{-a}^a  v^2
$$
and, using $v$ even in $x$ and in $y$ 
$$ E(v)= E_x(v) + E_y(v), \qquad
E_x(v) = 2\int_{-a}^a v(x,b)^2\ dx,\qquad
E_y(v) = 2\int_{-b}^b v(a,y)^2\ dy .
$$

When 
\begin{equation}
v_r(x,y)=\xi(r x)\eta(y/r)
\label{eq:vrdef}
\end{equation}
we write $F(v)=f_x(r) f_y(r)$ where
$$ f_x(r,a)= \int_{-a}^a \xi(r x)^2\, dx, \qquad
f_y(r,b)= \int_{-b}^b \eta(y/r)^2 \, dy .
$$
We have
$$ f_x(r,h) =\frac{1}{r}f_x(1,rh), \qquad
f_y(r,h) = r f_y(1,\frac{h}{r}) ,$$
and, in particular
$$ \int_{-h}^h  \int_{-h}^h \xi(rx)^2 \eta(y/r)^2 \, dx\, dy
=  \int_{-h/r}^{h/r}  \int_{-r h}^{r h} \xi(x)^2 \eta(y)^2 \, dx\, dy .
$$
Also define
$$ g_x(r,a)= \int_{-a}^a \xi'(r x)^2\, dx, \qquad
g_y(r,b)= \int_{-b}^b \eta'(y/r)^2 \, dy .
$$
We have
$$(\frac{\partial v_r}{\partial y})^2\
=\frac{1}{r^2}\, \xi(rx)^2 \eta'(y/r)^2
$$
so that
$$  \int_{-b}^b\int_{-a}^a (\frac{\partial v_r}{\partial y})^2\, dx\, dy
= G_x(r,a,b)
= \frac{1}{r^2} f_x(r,a) g_y( r,b) , $$
and
$$ G_x(r,h,h) = \frac{1}{r^2} G_x(1,rh,h/r) .$$
Similarly
$$ G_y(r,h,h) = \frac{1}{r^2} G_y(1,rh,h/r) .$$
The edge terms $E$ similarly involve the  $f$:
$$E_x(r,a,b)=2\eta(b/r)^2 f_x(r,a), \qquad
E_y(r,a,b)=2\xi(ra)^2 f_y(r,b) .
$$

We will be comparing $v_r(x,y)$, with the structure given in equation~(\ref{eq:vrdef}),
over the square side $2h$ with $v_1(x,y)$ over the rectangle
$(-rh,rh)\times(-h/r,h/r)$.

We now apply these formulae with $v$ motivated by the solution for a rectangle.
In the special case
$$\xi(r x)= \cos( \mu_X r x), \qquad \eta(y/r)= \cos( \mu_Y y/r),
$$
$G_x$ and $G_y$, as well as the $E_x$, $E_y$ and $F$
can be written in terms of $f_x$ and $f_y$.
For general $\mu_X$ and $\mu_Y$ with the abbreviations
$$ t_X= \tan(\mu_X r h), \qquad t_Y=\tan(\mu_Y h/r),
$$ that the integrals evaluate as follows:
\begin{itemize}

\item With $\Omega=(-rh,rh)\times{(-h/r,h/r)}$ and $v=\cos(\mu_X x)\cos(\mu_Y y)$  we have
\begin{align*} 
f_x
&=  r h+ \frac{t_X}{\mu_X(1+t_X^2)} ,\qquad
&f_y
&= \frac{h}{r} + \frac{t_Y}{\mu_Y(1+t_Y^2)}
  , \\
E_x
&= \frac{2 f_x}{1+t_Y^2}, \qquad
&E_y
&= \frac{2 f_y}{1+t_X^2}
 , \\
 G_x
 &= f_x \mu_Y^2 \left( \frac{h}{r}-\frac{t_Y}{\mu_Y(1+t_Y^2)}\right), \qquad
&G_y
&=  f_y \mu_X^2 \left( r h-\frac{t_X}{\mu_X(1+t_X^2)}\right)
 ,\\
 {\cal V}_x
 &= \frac{G_x+E_x/\beta}{F} ,
 &  {\cal V}_y
 &=\frac{G_y+E_y/\beta}{F}  
 , \\
  &= \mu_Y^2  -  \frac{2 (\mu_Y t_Y -1/\beta)}{f_y(1+t_Y^2)}  , 
 &
 &=  \mu_X^2  -  \frac{2(\mu_X t_X -1/\beta)}{f_x(1+t_X^2)} .
\end{align*}
We remark that, with $t_Y$ depending on $\mu_Y$ as above,
$$ \frac{\partial{\cal V}_x}{\partial\mu_Y}
=\frac{2(t_Y\mu_Y-1/\beta)}{(\mu_Y f_y (1+t_Y^2))^2}\,
\left(-t_Y +\mu_Y\frac{h}{r} (1+2\mu_Y\frac{h}{r})(1+t_Y^2)\right) .
$$
This is consistent with $\mu_Y^2$ being the eigenvalue for the strip width $2h/r$, 
when $t_Y\mu_Y=1/\beta$.\\

\smallskip
With $\Omega=(-rh,rh)\times{(-h/r,h/r)}$ and $v=\cos(\mu x/r)\cos(\mu r y)$  we have,
with $t=\tan(\mu h)$,$f_1=h + t/(1+t^2)$,
\begin{align*} 
f_x
&= r f_1 , \qquad
&f_y
&= \frac{f_1}{r}
  , \\
E_x
&= \frac{2 f_x}{1+t^2}, \qquad
&E_y
&= \frac{2 f_y}{1+t^2}
 , \\
G_x
&=  f_x \mu^2 r\left( h-\frac{t}{\mu(1+t^2)}\right)
&G_y
&= \frac{f_y \mu^2}{r} \left( h-\frac{t}{\mu(1+t^2)}\right), \qquad
 ,\\
 {\cal V}_x
 &= \frac{G_x+E_x/\beta}{F} ,
 &  {\cal V}_y
 &=\frac{G_y+E_y/\beta}{F}  
 , \\
  &= r^2 \mu^2  -  \frac{2 r^2 (r\mu t -1/\beta)}{f_y(1+t^2)}  , 
 &
 &=(\frac{\mu}{r})^2  -  \frac{2(\mu t/r -1/\beta)}{r^2 f_x(1+t^2)} .
\end{align*} 
In particular when $r\mu_*\tan(\mu_* h)=1/\beta$, $\mu_X^2\le{r^2\mu_*^2}$.
When $\mu=\mu_\square$ so $\mu t=1/\beta$, we find
\begin{equation}
\lambda_1(r)\le{\cal V}_x +{\cal V}_y
=\lambda_\square\left(
 \frac{1}{2}\left( 
r^2+\frac{1}{r^2}\right) -
\frac{\beta(\sqrt{r}-\frac{1}{\sqrt{r}})^2 (r+1+\frac{1}{r})}
{h(1+\frac{\beta^2\lambda_\square}{2}) +\beta}
\right) .
\label{eq:VxVymuSq}
\end{equation}
This implies the weaker inequality
\begin{equation}
\lambda_1(r)
\le 
\frac{\lambda_\square}{2}\left( 
r^2+\frac{1}{r^2}\right) ,
\label{eq:VxVymuSqSimple}
\end{equation}
and we remark there is equality in this when $\beta=0$.

\item With $\Omega$ as the square $(-h,h)\times{(-h,h)}$ and  $v=\cos(\mu_X r x)\cos(\mu_Y y/r)$ 
the integral defining $F$ evaluates to the same value as 
for the rectangle above with the same $v$ and
$$f_{x\square}= \frac{1}{r}  f_x , \qquad f_{y\square}= r f_y . $$
The other integrals, subscripted here with a $\square$ relate to those above as
\begin{align*}
E_{x\square}
&=  \frac{2 f_{x\square}}{1+t_Y^2}=\frac{1}{r} E_x ,\qquad
&E_{y\square}
&= \frac{2 f_{y\square}}{1+t_X^2} ={r} E_y 
  , \\
 G_{x\square}
  &=\frac{ f_{x\square} \mu_Y^2}{r} \left( \frac{h}{r}-\frac{t_Y}{\mu_Y(1+t_Y^2)}\right) 
&G_{y\square}
  &= r  f_{y\square} \mu_X^2 \left( r h-\frac{t_X}{\mu_X(1+t_X^2)}\right)\\
&= \frac{1}{r^2} G_x ,\qquad
&
&= {r^2} G_y 
,\\
 {\cal V}_{x\square}
 &= \frac{G_{x\square}+E_{x\square}/\beta}{F} ,
 &  {\cal V}_y
 &=\frac{G_y+E_y/\beta}{F}  
 , \\
  &=(\frac{ \mu_Y}{r})^2  -  \frac{2 ((\frac{ \mu_Y}{r}) t_Y -1/\beta)}{f_{y\square}(1+t_Y^2)}  , 
 &
 &= (r \mu_X)^2  -  \frac{2(r\mu_X t_X -1/\beta)}{f_{x\square}(1+t_X^2)} . 
\end{align*} 
Once again we can differentiate ${\cal V}_{x\square}$ with respect to $\mu_Y$. 
This time the derivative is zero when $\mu_Y t_Y/r =1/\beta$, consistent with what we know about
the fundamental eigenvalue for a strip of width $h$, as previously denoted $\mu_\square^2$.
\end{itemize} 

When $\mu_X$ and $\mu_Y$ satisfy the transcendental equations~(\ref{eq:muXmuY})
the integrals evaluate  as follows.
\begin{itemize}

\item With $\Omega=(-rh,rh)\times{(-h/r,h/r)}$ and $v=\cos(\mu_X x)\cos(\mu_Y y)$  we have
\begin{align*}
f_x
&= \frac{{\beta}/{r}+h+\mu_X^2 h\beta^2}{1+\beta^2\mu_X^2} ,\qquad
&f_y
&= \frac{{\beta}{r}+h+\mu_Y^2 h\beta^2}{1+\beta^2\mu_Y^2}
  , \\
&=  h+ \frac{{\beta}/{r}}{1+\beta^2\mu_X^2} ,\qquad
 &
&= h + \frac{{\beta}{r}}{1+\beta^2\mu_Y^2}
  , \\
E_x
&= \frac{2 f_x\mu_Y^2  \beta^2 r}{1+\beta^2\mu_Y^2}, \qquad
&E_y
&= \frac{2 f_y\mu_X^2 \beta^2 /r}{1+\beta^2\mu_X^2}
 , \\
 G_x
 &= \frac{ f_x (-r\beta+h+\mu_Y^2 h\beta^2)\mu_Y^2}{1+\beta^2\mu_Y^2}, \qquad
&G_y
&= \frac{ f_y (-{\beta}/{r}+h+\mu_X^2 h \beta^2)\mu_X^2}{1+\beta^2\mu_X^2}
 ,\\
&={ f_x \mu_Y^2}\left( h- \frac{r\beta}{1+\beta^2\mu_Y^2}\right), \qquad
&
&={ f_y \mu_X^2}\left( h- \frac{r/\beta}{1+\beta^2\mu_X^2}\right)
 ,\\
 \mu_Y^2
 &= \frac{G_x+E_x/\beta}{F}, \qquad
 &\mu_X^2
 &= \frac{G_y+E_y/\beta}{F} .
\end{align*}
This is consistent with 
$$\lambda_1({\rm rectangle}) = {\cal V}(v) =  \mu_X^2 + \mu_Y^2 .$$

\item With $\Omega$ as the square $(-h,h)\times{(-h,h)}$ and  $v=\cos(\mu_X r x)\cos(\mu_Y y/r)$ where
$\mu_X$ and $\mu_Y$ are as above, the integrals defining $f_x$ and $f_y$ evaluate to the same values as above.
The other integrals, subscripted here with a $\square$ relate to those above as
\begin{align*}
E_{x\square}
&= \frac{1}{r} E_x ,\qquad
&E_{y\square}
&= {r} E_y 
  , \\
 G_{x\square}
&= \frac{1}{r^2} G_x ,\qquad
&G_{y\square}
&= {r^2} G_y .  
\end{align*}
The variational characterization of $\lambda_{1\square}$ for a square gives
\begin{equation}
 2\mu_\square^2
= \lambda_{1\square}
\le {\cal V}_\square(v)
= \frac{G_x/r^2+E_x/(r\beta)}{F} + \frac{G_y/r^2+E_y/(r\beta)}{F} .
\label{eq:ineqVSqR}
\end{equation}
\end{itemize}

We believe that it may be possible to give a proof of Theorem~\ref{thm:thm3} beginning from
the inequality of~(\ref{eq:ineqVSqR}).
So far we have only numerical evidence that this is the case.
Indeed it may even be that this variational approach improves on the result based on $\phi_1$ that we
found at the end of the preceding subsection, and if so, we would have:
$$2\mu_\square^2 \le {\cal V}(v) \le \frac{1}{2}(\mu_X+\mu_Y)^2 
\le \frac{1}{2}\left((\mu_X+\mu_Y)^2+(\mu_X-\mu_Y)^2\right)= \mu_X^2 + \mu_Y^2 .$$
At present we do not have any useful upper bound on ${\cal V}(v)$.

\subsection{Inequalities on $\tan$ and an approach to Theorem~\ref{thm:thm3} in which $\beta$ and $h$ are specified.}
\label{subsec:hrSpecified}

\begin{figure}[hb]
\centerline{
$$
\includegraphics[height=5cm,width=6.5cm]{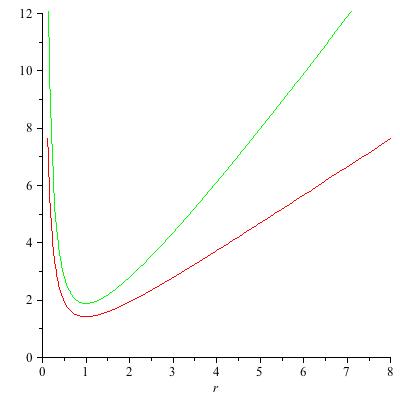}
$$
}
\caption{$\beta=1$ and $h=1$. For rectangles with area 4, $\lambda_{\rm LB}$ and
$\lambda_{\rm UB}$ are minimized by the square.
The horizontal axis is the $r$-axis.}
\label{fig:lambda1rectLBUB}
\end{figure}

With the definitions of equations~(\ref{eq:LBUBmu}),
define $\lambda_{\rm LB}(r)=\mu_{\rm LB}( h r)^2 +\mu_{\rm LB}( h/r)^2$ and
$\lambda_{\rm UB}(r)=\mu_{\rm UB}( h r)^2 +\mu_{\rm UB}( h/r)^2$.
Calculations show that both $\lambda_{\rm LB}(r)$ and $\lambda_{\rm UB}(r)$
are minimized at $r=1$. 
A typical plot is shown in Figure~\ref{fig:lambda1rectLBUB}.
Calculations show that both $\lambda_{\rm LB}(r)$ and $\lambda_{\rm UB}(r)$
are decreasing functions of $r$ for $0<r<1$ and increasing for $r>1$.
The calculations for $\lambda_{\rm LB}(r)$ are the easier, and one finds
$\lambda_{\rm LB}(r)$ is convex on $(0,\infty)$, but not log-convex.

\smallskip

We return now to how this might lead to an alternative proof of Theorem~\ref{thm:thm3}. 
The present strategy in this approach is to consider, separately, ranges away from $r=1$ (Result 3A) and
a range about $r=1$ (Result 3B).

\par\noindent{\bf Result 3A.}
{\it Amongst all rectangles with a given area, those which are long and thin have larger fundamental Robin eigenvalues
than the square, and there is an explicit bound on the aspect ratios (denoted $r_*^2$ below and its reciprocal)
which is given below.}

 \smallskip
 
 \par\noindent{\it Comment.} 
 As $\lambda_1$ tends to infinity as $r$ tends to 0 or as $r$ tends to infinity,
 the first part of the last sentence of the theorem statement is trivial..
 The nontrivial part is the bound the intervals on which we have shown the result.

\medskip
\par\noindent{\it Proof of Result 3A.} 
From the formulae for $\mu_{\rm LB}(c,\beta)$ and $\mu_{\rm LB}(c,\beta)$ given in
equations~(\ref{eq:LBUBmu}), and especially their consequences for the properties of $\lambda_{\rm LB}$,
 we can deduce that $\lambda_1(r)\ge\lambda_1(1)$ for the range of $\beta$,
intervals $(0,r_*)$ and $(1/r_*,\infty)$  where
$${\frac {{\pi }^{2} \left( 4\,h+{\pi }^{2}\beta\,r+4\,{r}^{4}h+{r}^{3}{
\pi }^{2}\beta \right) }{hr \left( 4\,hr+{\pi }^{2}\beta \right) 
 \left( 4\,h+{\pi }^{2}\beta\,r \right) }}
=\lambda_{\rm LB}(r) >  \lambda_{\rm UB}(1) .$$ 
The expression for $\lambda_{\rm LB}(r)$ on the left is the ratio of a quadratic function of
$(r+1/r)$ in the numerator to a linear function of $(r+1/r)$ in the denominator
(as noted in equation~(\ref{eq:qrPlus1or})).
Finding $r_*(\beta/h)<1$ involves solving a quartic $q(r,y)=0$,
and we begin with the more general question of finding the $r$ where
$\lambda_{\rm LB}(r)=y/h^2$, for 
\begin{equation}
\frac{y}{h^2}>\lambda_{\rm LB}(1)=  \frac{2\pi^2}{h(4 h +\pi^2 \beta)} . 
\label{eq:yCond}
\end{equation}
The quartic has the structure
\begin{equation}
 q(r,y) = r^4 + a_1 r^3 + a_2 ^2 + a_1 r + 1 =0 ,
\label{eq:quarticar}
\end{equation}
and the dependence on $y$ is linear: $q(r,y)=q_0(r)-y q_1(r)/h^2$ with $q_0(r)=q(r,0)$.
The quartic has the property $q(1/r,y)=q(r,y)/r^4$, i.e. if $r=R$ is a root, so is $r=1/R$.
As
\begin{equation}
 \frac{q(r,y)}{r^2}
= \left( r+\frac{1}{r}\right)^2 + a_1 \left( r+\frac{1}{r}\right) + (a_2-2) .
\label{eq:qrPlus1or}
\end{equation}
it is straightforward to first solve for $\left( r+\frac{1}{r}\right)$, and then for $r$.
For the actual quartic in the application
\begin{equation}
 a_1 =  \frac{\beta \pi^2}{4 h} - y \frac{\beta}{h},\qquad
a_2= -y\, \frac{\beta^2\ \pi^4 + 16 h^2}{4 h^2 \pi^2} .
\label{eq:a1a2def}
\end{equation}
The equation $\lambda_{\rm LB}(r)=y/h^2$ is
$$ \frac{q_0(r)}{q_1(r)} = \frac{y}{h^2}. $$

The quadratic in $\left( r+\frac{1}{r}\right)$ given in
equation~(\ref{eq:qrPlus1or}) has zeros at
$$\left( r+\frac{1}{r}\right)
= \frac{1}{2}
\left( -a_1 \pm \sqrt{a_1^2 -4 a_2 +8}\right) .
$$
From equation~(\ref{eq:a1a2def}), $a_2<0$ so that only the positive root above is
appropriate: denote it by $B$. Then the solutions for $r$ are
$$ r_\pm = \frac{1}{2}\left( B \pm \sqrt{B^2 -4}\right) . $$
We remark that $B>2$ (and obviously $\left( r+\frac{1}{r}\right)>2$ on $(0,\infty))$.
That $B>2$ is equivalent to $-4 a_2> 8 +2 a_1$ which in turn is
equivalent to inequality~(\ref{eq:yCond}).

This completes the proof of Result 3A.

\medskip
\par\noindent{\bf  Range about $r=1$.}

\par\noindent{\bf Result 3B.}
{\it  For values of $r$ close to 1, $\lambda_1(r)\ge\lambda_1(1)$.}

\par\noindent{\it Proof.}
We have set up coupled first order differential equations for $\lambda(r)$ and $\delta(r)$
where
$$ \lambda(r)=\mu_{(2)}(r)+\mu_{(2)}(\frac{1}{r}), \qquad
\delta(r)=\mu_{(2)}(r)-\mu_{(2)}(\frac{1}{r}) .
$$
We already know from properties of $\mu$ that $\delta(r)$ is a decreasing function for $r\in(0,\infty)$.
Series expansions about $r=1$ can be found and locally 
$$ \lambda(r)\sim \lambda(1) + L_2 (r-1)^2, \qquad {\mbox{\rm with\ }} L_2>0 \qquad
{\mbox{\rm as\ }}r\rightarrow{1} .
$$
This shows that, very locally about $r=1$, $\lambda_1(r)>\lambda_1(1)$.
With further work it is possible to determine an interval.
\medskip
\par\noindent
We believe it would be possible to give an alternative proof of Theorem~\ref{thm:thm3} by establishing
intervals in the two ranges above, 3A and 3B, so that they overlap.
It may be that other upper or lower bounds than those used in this subsection
(e.g. some as suggested in~\S\ref{subsec:variational}) may be useful in this.

\clearpage

\section{Completely monotone functions and related topics}\label{sec:compmon}

Alternate proofs, in the notation of \S\ref{sec:rectlambda1} that $\lambda_1(r)\ge\lambda_1(1)$,
may be worth considering.
Some of these might depend on properties of $\mu(r)$.

This appendix has general facts concerning
completely monotone and absolutely monotonic functions and subsets thereof that
might be of relevance to our study of $\mu(c)$.

We wish to find alternative proofs of Theorem~\ref{thm:thm3}
\begin{equation}
 \mu_{(2)}(r) +  \mu_{(2)}(\frac{1}{r}) \ge 2  \mu_{(2)}(1) .
\label{eq:B1}
\end{equation}
Now
$$ \mu_{(2)}(r) +  \mu_{(2)}(\frac{1}{r}) \ge 2  \mu_{(2)}(\frac{1}{2}(r+\frac{1}{r}))) ,
\eqno{(\rm convexity)}$$
Rather more is true as $\mu$ and hence $\mu_{(2)}$ are log-convex:
$$ \sqrt{\mu_{(2)}(r) \   \mu_{(2)}(\frac{1}{r})} \ge 2  \mu_{(2)}(\frac{1}{2}(r+\frac{1}{r}))) 
\eqno{(\rm logconvexity)}$$

It may be that $\mu$ and hence $\mu_{(2)}$ are completely monotone.
We wish to see how complete monotonicity might allow the above two can be improved.
However we note below that inequality~(\ref{eq:B1}) is not satisfied by every completely monotone function.

\subsection{Log-convex functions}

\par\noindent
{\bf Lemma LC1.} {\it The function $f$ is log-convex on an interval $I$, if and only if for all 
$a, b, c \in{ I}$ with $a < b < c$, the following holds:
$$f(b)^{c-a} \le  f(a)^{c-b} \le f(c)^{b-a} . $$
}
\medskip
Hence for $r\le{1}\le{1/r}$,
$$ 1 \le \left( \frac{\mu(r)}{\mu(1)}\right)^{\frac{1}{r}-1} \left( \frac{\mu(\frac{1}{r})}{\mu(1)}\right)^{1-r} . $$

\vspace{1cm}

\subsection{Completely monotone functions}

\par\noindent {\bf Definition.}
A real-valued function $f$ defined on $[0,\infty)$ is said to be
{\it completely monotone} (totally monotone, completely monotonic, totally monotonic)
  if   
 $(-1)^k f^{(k)}(x)\ge{0}$ for $x>0$ and $k=0,1,2,\ldots$.\\
Denote the set of completely monotonic functions on $[0,\infty)$ by ${\CM}$.
\medskip

\par\noindent
{\bf Theorem CM1.} {\it
The set ${\CM}$ forms a convex cone: $(t_1 f_1 + t_2 f_2)\in{\CM}$
for all nonnegative numbers $s$, $t$ and all $f_1\in{\CM}$ and $f_1\in{\CM}$.\\
The set ${\CM}$ is also closed under multiplication and point-wise convergence. That is
$$f_1(x) f_2(x) \in {\CM}\qquad{\rm	and\ \ }
\lim_{n\rightarrow\infty} f_n(x) \in{\CM}, $$
where $f_n(x)\in{\CM}$ for all $n\ge{1}$ and their point-wise limit exists for any $x > 0$.}
\medskip

\par\noindent
{\bf Theorem CM2.} {\it Let $f(x)\in{\CM}$. and let $h(x)$ be nonnegative with its derivative in ${\CM}$.
Then $f(h(x)\in{\CM}$.}
\medskip

\par\noindent
{\bf Corollary CM2.} {\it Let $f(x)\in{\CM}$ and $f(0) < \infty$. Then the function
$$ -\log\left( 1=\frac{f(x)}{A}\right) . \qquad
A\ge{f(0)}, 
$$
is ${\CM}$. From this it follows that
$$\frac{f'(x)}{A-f(x)}. \qquad A\ge f(0) $$
is ${\CM}$ since this reduces to minus the derivative of the previous expression.
}
\smallskip

\par\noindent This corollary is given in~\cite{MS01}.

\medskip

From the derivative of the last function, we have
$$\frac{d}{d x}  \frac{f'(x)}{A-f(x)}
= \frac{ f'' (A-f) + (f')^2}{(A-f)^2} \le {0} . $$
Rearranging this gives
$$ f f'' - (f')^2\ge A f'' \ge{0} . $$
In particular, any $f\in\CM$ is log-convex.
\smallskip

There is a relationship of $\CM$ functions and Laplace transforms.
Define
$$F(r)= \int_0^\infty \exp(-r t) f(t) \, dt .$$
If the above integral is bounded for all $r>0$ and $f(t)\ge{0}$ for all $t>0$
then $F\in\CM$.
The converse is also true.
\smallskip

That the Laplace transform of a positive function is log-convex can be proved using the
Cauchy-Schwarz inequality.  Suppose $f(t)\ge{0}$. Then
\begin{eqnarray*}
F(\frac{r_1+r_2}{2})
&=&\int_0^\infty \exp(-\frac{r_1+r_2}{2}\, t) f(t)\, dt \\
&=&  \int_0^\infty \exp(-\frac{r_1}{2}\, t) \sqrt{f(t)}\,  \exp(-\frac{r_2}{2}\, t) \sqrt{f(t)} \, dt
\le \sqrt{F(r_1)\, F(r_2)} .
\end{eqnarray*}
See \S4.8 of~\cite{Wa92}.

\medskip

If we were to attempt to establish (\ref{eq:B1}) for the Laplace transforms of positive functions, one
might begin with
$$ F(r)+F(\frac{1}{r})-2 F(1)
= \int_0^\infty k(r,t) f(t)\, dt ,
$$
where
$$ k(r,t)=\exp(-r t) + \exp(-\frac{t}{r}) - 2\exp(-t) . $$
Also define
$$ k_c(r,t)=\exp(-r t) + \exp(-\frac{t}{r}) - 2\exp(-\frac{1}{2}(r+\frac{1}{r})\, t) \ (\ge{k(r,t)})\ . $$
Now $k_c(r,t)\ge{0}$ for $r>0$ and $t>0$ and this is another way to show $F$ is convex.
It happens that, for every $r>0$,  $k(r,t)$ takes on both signs, and, as a consequence, there are functions $F\in\CM$
which do not satisfy~(\ref{eq:B1}), so, if $\mu_{(2)}$ does it is as a consequence of further properties.
At fixed $r>0$, $k(r,t)<0$ for $0<t<1/2$ and $k(r,t)>0$ for $t>1$, and there is a unique $t_{0}(r)$ in
the interval $(1/2,1)$ where $k(r,t_{0}(r))=0$.

\medskip

Another published statement, a special case of which is that
any $f\in\CM$ is log-convex is as follows:\\
{\it Let $f\in{\CM}$ . Then
$$ (-1)^{nk}\left( f^{(k)}(x)\right)^n \le (-1)^{nk} \left( f^{(n)}(x)\right)^k (f(x))^{n-k} $$
for all $x>0$ and integers $n\ge{k}\ge{0}$.}\\
In particular, for $n=2$ and $k=1$ we have that any completely monotonic function is log-convex.

\subsection{Absolutely monotonic  functions}

A function $f(x)$ is {\it absolutely monotonic} in the interval $a<x<b$ if it has nonnegative derivatives of all orders in the region, i.e.,
$f^{(k)}(x)\ge{0}$ for all $x$ in the interval and $k=0,1,2,\ldots$\ . 
Denote by $\AM{(a,b)}$ the set of all functions {absolutely monotonic} in the interval $a<x<b$.
\medskip

\par\noindent
{\bf Theorem AM1.} {\it
The set ${\AM}{(a,b)}$ forms a convex cone: $(t_1 f_1 + t_2 f_2)\in{\AM}{(a,b)}$
for all nonnegative numbers $s$, $t$ and all $f_1\in{\AM}$ and $f_2\in{\AM}{(a,b)}$.\\
The set ${\AM}$ is also closed under multiplication and point-wise convergence. That is
$$f_1(x) f_2(x) \in {\AM}{(a,b)}\qquad{\rm	and\ \ }
\lim_{n\rightarrow\infty} f_n(x) \in{\AM}{(a,b)}, $$
where $f_n(x)\in{\AM}{(a,b)}$ for all $n\ge{1}$ and their point-wise limit exists for any $x > 0$.}
\medskip

In particular the function $X\tan(X)$ is absolutely monotonic on $[0,\pi/2)$.
\medskip

Widder (1941)~\cite{Wi41} gives:

\par\noindent
{\bf Theorem AM2.} {\it
$f\in\AM$ and $g\in\CM$ then the composition $f\circ{g}\in\CM$.
}

\subsection{Miscellaneous topics}

Further classes of functions, Stieltjes functions, Bernstein functions, etc., are treated in~\cite{SSV}.

Let {\bf CMI} tbe the set of functions $\phi$ with domain and range $(0,\infty)$, 
with  $\phi\in\CM$, and for which the inverse $\phi^{-1}$ is also in $\CM$.
The set {\bf CMI} is nonempty as
$f(\alpha,z)=z^{-\alpha}$ is in $\CM$ when $\alpha>0$ and
$f(1/2,f(2,z))=z$ so $f(1/2,\cdot)$ has as its inverse $f(2,\cdot)$.
The question arises as to whether either of the functions $\phi_1$ or $\phi_2$ defined in \S\ref{subsec:Proofexplicit} is in {\bf CMI}.
Evidence that this might be the case is noted in \S\ref{subsec:Rectangular}.
There are two questions related to this.\\
(i) If the $\phi$ are in {\bf CMI} how might this help with alternative proofs of Theorem~\ref{thm:thm3}?\\
(ii) Under what conditions is the inverse of a function $f\in\CM$ also in $\CM$? And, do our $\phi$ satisfy these conditions?

\clearpage

\end{document}